\documentclass[authoryear,preprint,11pt]{elsarticle}
\usepackage{multirow}
\usepackage{pgf-pie}
\usepackage{color,xcolor}
\usepackage[utf8]{inputenc}
\usepackage{tikz}
\usepackage{enumitem}
\usepackage{amsmath,amsfonts,amsthm}

\usepackage{subcaption}

\usepackage{gensymb}
\usetikzlibrary{shapes.geometric, arrows, positioning, fit, shapes, arrows.meta}
\tikzstyle{bigbox} = [draw, dashed, rounded corners, text width=18em, minimum height=8em, align=right]

\tikzstyle{legend0} = [draw, rectangle, fill=red!5, text width=0.55cm, align=center]
\tikzstyle{legend1} = [draw, rectangle, fill=red!10, text width=0.55cm, align=center]
\tikzstyle{legend2} = [draw, rectangle, fill=red!20, text width=0.55cm, align=center]
\tikzstyle{legend3} = [draw, rectangle, fill=red!30, text width=0.55cm, align=center]
\tikzstyle{legend4} = [draw, rectangle, fill=red!40, text width=0.55cm, align=center]
\tikzstyle{legend5} = [draw, rectangle, fill=red!50, text width=0.55cm, align=center]
\tikzstyle{legend6} = [draw, rectangle, fill=red!60, text width=0.55cm, align=center]
\tikzstyle{legend7} = [draw, rectangle, fill=red!70, text width=0.55cm, align=center]
\tikzstyle{legend8} = [draw, rectangle, fill=red!80, text width=0.55cm, align=center]

\tikzstyle{bigbox1} = [draw, dashed, rounded corners, text width=18em, minimum height=4em, align=right]
\tikzstyle{startend} = [rectangle, draw, fill=!40, text width=3em, text centered, rounded corners, minimum height=2em]
\tikzstyle{block} = [circle, draw, fill=gray!50, text width=7.5em, text centered, minimum height=2em, font]
\tikzstyle{block2} = [circle, draw, fill=white!25, text width=7.5em, text centered, minimum height=2em, font]
\tikzstyle{block3} = [rectangle, text centered, rounded corners]
\tikzstyle{block6} = [rectangle, text width=40em, text centered, rounded corners, minimum height=2em]
\tikzstyle{block4} = [rectangle, draw, fill=gray!15, text width=6em, text centered, rounded corners, minimum height=2em]
\tikzstyle{block5} = [rectangle, draw, text width=10em, text centered, minimum height=2em]
\tikzstyle{blockns} = [rectangle, draw, text width=5em, text centered, minimum height=5em]
\tikzstyle{blockns1} = [rectangle, draw, text width=2em, text centered, minimum height=2em]
\tikzstyle{blocks} = [rectangle, draw, fill=gray!50, text width=5em, text centered, minimum height=5em]
\tikzstyle{blocks1} = [rectangle, draw, fill=darkgray!15, text width=2em, text centered, minimum height=2em]
\tikzstyle{block0} = [rectangle, draw, fill=gray!22, text width=5em, text centered, minimum height=5em]
\tikzstyle{decision} = [diamond, text width=2cm, minimum width=3cm, minimum height=1cm, text centered, draw, fill=white]
\tikzstyle{line} = [draw, -latex']
\usepackage{makecell}
\usepackage{amssymb}
\usepackage{eurosym}
\usepackage{color}
\usepackage{bm}
\usepackage{booktabs}
\usepackage{multirow}
\usepackage{longtable}
\usepackage{arydshln}
\usepackage[left = 1.7cm, right = 1.7cm, top = 1.7cm, bottom = 1.7cm]{geometry}
\usepackage{enumitem}
\usepackage{algorithm} 
\usepackage{algpseudocode} 
\usepackage{subcaption}
\usepackage{tikz}
\usepackage{pgfplots}
\usetikzlibrary{patterns}
\usetikzlibrary{arrows}
\usepackage{lscape}

\usepackage{amsthm,amsmath,eurosym}

\journal{European Journal of Operational Research}

\usepackage{lineno, blindtext}

\makeatletter
\def\ps@pprintTitle{%
  \let\@oddhead\@empty
  \let\@evenhead\@empty
  \def\@oddfoot{\reset@font\hfil\thepage\hfil}
  \let\@evenfoot\@oddfoot
}
\makeatother

\begin{document}

\begin{frontmatter}
    \title{New Multi-objective Partial Optimisation Decomposition Strategies\\ for the Thesis Defence Scheduling Problem}
    \author[cegist]{João {\sc  Almeida}\corref{cor1}}\ead{joao.carvalho.almeida@tecnico.ulisboa.pt}

    \author[inesc]{Alexandre P. {\sc Francisco}}
    
    \author[cegist]{Daniel {\sc Santos}}
    
    \author[cegist]{Jos\'e Rui {\sc  Figueira}}

    \address[cegist]{CEGIST, Instituto Superior T\'{e}cnico,  Universidade de Lisboa, Portugal}

     \address[inesc]{INESC-ID, Instituto Superior T\'{e}cnico,  Universidade de Lisboa, Portugal}

    \cortext[cor1]{Corresponding author at: CEGIST, Instituto Superior T\'{e}cnico,  Universidade de Lisboa, Portugal.}

    \begin{abstract}
    \noindent A new multi-objective method for the thesis defence scheduling problem is introduced. This is an academic scheduling problem with a real-world impact on the operations of universities, the daily lives of committee members, and the quality of the assessment process. The problem consists of appointing committees to defences and assigning them to a time slot and room. A multi-objective approach is necessary to provide a better understanding of the possible solutions and their trade-offs to the decision-makers. However, this type of approach is often time-consuming. The new multi-objective optimisation approach decomposes the monolithic problem into a sequence of multi-objective problems. This leads to significant efficiency gains when compared to the augmented-$\epsilon$ constraint method. In our application, the monolithic model is decomposed into two submodels (stages) to be solved sequentially. In the first stage, genetic algorithms (NSGA-II and NSGA-III) are used to find multiple committee configurations (partial solutions). The performance of these solutions is assessed based on committee assignment quality objectives and a proxy objective which aims to predict the performance of the objectives in the next stage. In the second stage, considering multiple partial solutions found previously, an augmented $\epsilon$-constraint method is solved to find non-dominated solutions regarding the assignment of time slots to each defence. These solutions consider schedule quality objectives. Finally, non-dominated solutions are presented based on objective function performance for both points of view. For small size instances, the method takes between 8\% and 32\% of the time of an augmented $\epsilon$-constraint method but finds non-dominated sets with slightly worse hyper-volume indicator values. For larger size instances, the times are between 6\% and 18\% of monolithic resolutions and the hyper-volume indicator values are better. A real-world case study is presented. The experiment with decomposition found 39 non-dominated solutions in 1600 seconds. The augmented $\epsilon$-constraint method found 9 solutions in 2400 seconds. For the three objectives, the new method found a solution that improved the best-performing solution found with the other method in the time limit being considered.

    \end{abstract}
    \vspace{0.25cm}
\end{frontmatter}




\section{Introduction}\label{sec:introduction}

\noindent  Organising thesis defences where university students present and defend their work before a committee of experts, can be a challenging operational task for universities and their departments. The problem involves appointing committee members to assess a defence and then scheduling defences to a time slot considering the availability of the committee members and rooms. Doing this process without optimisation tools is inefficient and can generate unfair timetabling arrangements \citep{BattistuttaEtAl2019, Almeida2024}.  The problem can be characterised as assigning committee members to perform predetermined roles in a set of thesis defences/committees. An appropriate time slot and room must also be assigned for each defence.

The people responsible for conducting this process at the university level struggle to find feasible assignments. However, when they are possible, their quality is often underwhelming, leading to dissatisfaction among the committee members. From the point of view of the committee assignment, there can be some problems regarding the unfair distribution of workloads or a bad match between the expertise of the committee members and the subjects of the defences. From the point of view of the schedules, problems such as an excessive amount of days with scheduled defences can significantly increase the impact this responsibility has on the daily lives of the committee members.

When dealing with problems with a single objective function, the methods aim to find one solution that has the best performance considering those functions. Multi-objective optimisation problems arise if there are conflicting objectives to be considered when a certain decision is made. The aim is to find multiple solutions with trade-offs between different objective function performances. This provides decision-makers with a better knowledge of the solution space and allows better-informed decisions. For example, \cite{Zeitrag2024} deals with a multi-objective lot sizing and production scheduling of a real-world automotive industry problem. Their approach provides decision-makers with non-dominated solutions regarding workload balance, total cost and average inventory.

There are heuristic multi-objective methods, such as genetic algorithms, and deterministic methods, such as the $\epsilon$-constraint method. The principles of natural selection inspire genetic algorithms. A population of candidate solutions to a problem is generated and iteratively evolves through selection, crossover, and mutation operations to produce better solutions over successive generations. They are applied to several problems in the literature, such as job shop scheduling problems \citep{Zeitrag2022, Tutumlu2023}, travelling salesman problems \citep{Zheng2023, Mahmoudinazlou2024, Mahmoudinazlou20241, Soares2024}, vehicle routing problems \citep{Bortfeldt2020, Xue2021, Zhao2024}, or project scheduling problems \citep{Servranckx2024, Bredael2024}. The $\epsilon$-constraint method iteratively solves a single objective problem but sets lower or upper bounds for the performances of other objectives. This iterative procedure makes finding non-dominated solutions with trade-offs between their objectives possible. The method is used to solve several problems, such as task scheduling problems \citep{Stewart2023, Wu2024}, knapsack problems \citep{Mesquita-CunhaEtAl2022}, or the generalized cubic cell formation problem \citep{Bouaziz2023}.

To provide this additional information to decision-makers, multi-objective approaches are often more time-consuming than single-objective ones. This work proposes a new multi-objective partial optimisation method, combining elements of the $\epsilon$-constraint method and genetic algorithms. This new method is less time-consuming than other multi-objective optimisation techniques. 

Other multi-objective approaches typically optimise the entire problem simultaneously using a single model, known as a monolithic model. We propose a decomposition approach for multi-objective optimisation to enhance efficiency while maintaining high-quality solutions. Instead of considering a monolithic model, the problem is formulated as a sequence of submodels, or stages. The non-dominated solutions identified in a stage are used as partial solutions for further optimisation in the next stage. Essentially, each solution obtained in a prior stage serves as a basis for applying a new multi-objective algorithm to identify additional non-dominated solutions derived from each partial solution. To ensure that the objectives of the next stage are represented, proxy objectives are included in the previous stage. 

The thesis defence scheduling problem is a good fit to test this method as it has multiple conflicting objectives, such as maximising the fit between the expertise of the committees and the subject of the defences or minimising the number of days a committee member is scheduled for. Moreover, it can be decomposed into two stages, the committee assignment and the time-slot appointment. During the first stage, a genetic algorithm (NSGA-II or NSGA-III) is employed to find a set of committee configurations, serving as partial solutions for the subsequent stage. In the following stage, an augmented $\epsilon$-constraint method is applied to each configuration to determine non-dominated time-slot assignments for every defence. All solutions generated for each committee are collectively evaluated to identify those that remain non-dominated across the entire set, and these are then presented to the decision-makers.

The main novelty in this pa\textit{per} is the new method and its application to a problem with real-world situations. The method differs from other decomposition methods due to its multi-objective nature and the inclusion of proxy objectives in the first stages. Additionally, it is distinct from other multi-objective approaches due to its multi-stage structure, which improves its efficiency.

The method is tested in small and large size randomly generated instances and a real-world case study. For the smaller size instances, where the augmented $\epsilon$-constraint iterations could reach optimality within a 120 seconds time limit, the decomposed method takes between 8\% and 32\% of the time of monolithic resolutions, but the monolithic resolutions outperform it regarding the quality of the solutions. This is expected as this is a comparison between a partial optimisation technique and a deterministic method which can reach optimality. For the larger size instances, where the augmented $\epsilon$-constraint iterations could not reach optimality, the decomposed method takes between 6\% and 18\% of the time of monolithic resolutions and is capable of finding non-dominated sets with larger hyper-volume indicator values. The real-world case study illustrates the performance of three objectives, one which is evaluated in the first stage and two being evaluated in the second stage. The first stage objective and one of the second stage objectives have very similar performances when compared to a monolithic resolution. The other second stage objective had a slightly worse performance. This is due to one of the second stage objectives having a stronger relationship with the proxy objective than the other one.
 
The remainder of this pa\textit{per} is organised as follows. Section \ref{s-lr} provides a literature review regarding decomposition and partial optimisation methods and the thesis defence scheduling problem. Section \ref{s-prob} gives an introduction and motivation to study thesis defence scheduling. Sections \ref{s-mono} and \ref{s-dec} present the monolithic and decomposed models, respectively. Section \ref{s-augm} introduces some theoretical background regarding $\epsilon$-constraint methods. Section \ref{s-gen} explains the genetic algorithms used in our application. Section \ref{s-AlgF} provides an overview of the algorithmic framework. Sections \ref{s-exp} and \ref{s-cases} discuss the results considering the randomly generated instances and the real-world case study, respectively. Section \ref{conc} concludes the work and proposes future research avenues.

\section{Literature review}\label{s-lr}

\noindent This section presents a literature review on the scope of our work. It presents works regarding decomposition and partial optimisation and the thesis defence scheduling problem.

\subsection{Decomposition and partial optimisation}\label{dpo}

\noindent \cite{BurkeEtAl2010Decomposition} define the decomposition method for integer programming problems. The pa\textit{per} proposes an approach to solving optimisation problems with multiple interacting components, each associated with different sets of objectives. It can be defined as the sequential optimisation of multiple restricted submodels. Initially, only one computationally difficult component and its associated subset of objectives are considered, leading to partial solutions that define interesting neighbourhoods in the complete problem's search space. Variable aggregation is performed in the first stage, allowing the exploration of neighbourhoods guaranteed to contain feasible solutions. The approach employs integer programming to implement heuristics that produce solutions with bounds on their quality. This work applies the method to the ITC2007 curriculum-based course timetabling problem. The approach has mixed results, which vary for different instances. In general, the state-of-the-art local search approaches at the time of the publication of this work outperformed it. Nonetheless, decomposing the problem improved the performance when compared to monolithic resolutions using other linear programming based methods.

\cite{Dunke2023} applied decomposition to the course timetabling problem of their university. This problem is decomposed into three submodels. In the first submodel, the time slots in which the lectures of certain study plans are defined. The second submodel defines the time slots for tutorial sessions. The last submodel assigns the individual student timetables. The work also proposes an initial step where a genetic algorithm enhanced with an artificial neural metamodel provides partial solutions which fix the number of lectures each study plan has each day. They have a good performance in two out of their three objective functions. Their method is efficient at solving smaller size instances but becomes exponentially more time-consuming for larger size instances as they solve their entire model for every partial solution found in the initial step.

\cite{Seizinger2023} addresses the shortage of nurses in industrialised countries, particularly in dual vocational training systems like those in Germany. It explores challenges in organising apprenticeship programs for nursing education and introduces two submodels to optimise planning levels. The study uses a real-world case from a German vocational school, emphasising the development of a method for strategic/tactical and tactical/operational planning. Their method can improve the supply of students for medical units by 40\% when compared to plans that are manually created.

\cite{Mikkelsen2022} decomposes the International Timetabling Competition 2019 course timetabling problem. This proved beneficial in some larger size instances. They separate the student sectioning and timetable development allowing them to solve instances that would be otherwise intractable for the solver.

\cite{Cataldo2017} addresses a real-world examination scheduling problem. The proposed solution involves four sequential submodels. The first clusters courses and generates room patterns, while the second assigns time slots and room patterns to course clusters. The third refines assignments for individual courses, and the fourth produces a definitive exam schedule with specific room assignments. The findings demonstrate a decrease in conflicts and rescheduling occurrences compared to the current exam scheduling practices employed at their university.

\cite{VermuytenEtAl2016} proposes a solution to a real-world course timetabling problem. In a first stage, they optimise the timetables disregarding room assignments. In a second stage, they assign the rooms with the aim of minimising congestion in the hallways as a consequence of student flows between classes. This approach produces near-optimal results in a much shorter time frame than the monolithic model. 

\cite{SorensenDahms2014} present a two-stage decomposition approach for addressing a real-world high school timetabling problem. The proposed decomposition divides the problem into submodels, with fewer variables. Irreversible decisions made in the first stage are influenced by the second stage objective, resembling a minimum weight maximum matching problem. The decomposition strategically utilizes Hall's theorem, requiring a manageable subset of constraints. The approach can find better solutions than a monolithic model in the time limit being considered.

\cite{Cacchiani2013} decomposes the International Timetabling Competition 2007 curriculum-based course timetabling problem. Different objectives are assessed in each of their submodels. This work proved that some of the heuristic solutions that had previously been found for this problem were optimal, and considerably outperformed other integer programming approaches.

One of the requirements of this method is that the problems have multiple components which translate into different objectives. However, to the best of our knowledge, previous works have aggregated the different objectives into weighted-sum objective functions and do not use this method to find non-dominated solution sets. In our work, we extend the decomposition approach to a multi-objective problem and use it to improve the efficiency of the multi-objective search.

\subsection{Thesis defence scheduling}\label{tds}

\noindent Three types of thesis defence scheduling problems have been identified in the literature \citep{Almeida2024}. The first type involves a single defence assignment, as addressed in works by \cite{HuynhEtAl2012}, \cite{PhamEtAl2015}, \cite{LimantoEtAl2019}, \cite{TawakkalSuyanto2020}, \cite{ChristopherWicaksana2021} and \cite{Almeida2024}. In these works, a committee is assigned to a single defence. The second type focuses on assignments regarding sessions of defences, exemplified by the works of \cite{BattistuttaEtAl2019} and \cite{Dimitsas2022}. In these works, each committee is assigned to a set of defences. The third type is the hybrid assignment type and is discussed in the study by \cite{KochanikovRudova2013}. This work involves a combination of the other scheduling approaches. Specifically, a committee is initially assigned to an extended period, similar to the session of defences assignment type. However, within this broader timeframe, individual thesis defences and any additional required committee members are assigned to specific hours.

To assess the scheduling scenarios, two primary points of view emerge, each associated with specific objectives \cite{Almeida2024}. Firstly, the committee assignment quality perspective incorporates objectives such as ensuring a fair distribution of assignments \citep{HuynhEtAl2012, KochanikovRudova2013, PhamEtAl2015} and optimising the expertise match between committee members and defence subjects \citep{HuynhEtAl2012, BattistuttaEtAl2019, ChristopherWicaksana2021}. Secondly, the schedule quality viewpoint encompasses objectives like promoting compact schedules \citep{HuynhEtAl2012}, satisfying preferred time slot requests \citep{TawakkalSuyanto2020}, and preventing room changes \citep{HuynhEtAl2012, BattistuttaEtAl2019}.

Several solution approaches have been explored to address these scheduling challenges: (1) Mixed-integer linear programming \citep{BattistuttaEtAl2019, Almeida2024}; (2) Constraint programming \citep{BattistuttaEtAl2019}; (3) The greedy backtracking hybrid algorithm \citep{TawakkalSuyanto2020}; (4) Constructive heuristic \citep{KochanikovRudova2013}; (5) Local search methods \citep{PhamEtAl2015, BattistuttaEtAl2019}; (6) Genetic algorithms \citep{HuynhEtAl2012, LimantoEtAl2019}; (7) Particle swarm optimisation \citep{ChristopherWicaksana2021}.

This work addresses the thesis defence scheduling problem presented in \cite{Almeida2024}, which can be applied to solve other single assignment thesis defence scheduling problems from the literature. Some improvements to the formulation have been made. Moreover, a new multi-objective solution approach is proposed. Compared with the augmented $\epsilon$-constraint implemented by \cite{Almeida2024}, the multi-objective decomposition approach is significantly more time efficient.

\section{The Thesis Defence Scheduling Problem}\label{s-prob}

\noindent To successfully conclude their academic degrees, most university students are required to present and defend their thesis before a committee of experts in the field. Assigning these committees and scheduling the defences by hand is a burdensome task and often leads to poor or unfair member timetable arrangements. 

When doing the assignment of committee members to defences, their fields of expertise should be matched to the subjects of the theses. Moreover, it is also important to be fair when considering the workloads of the committee members and attempting to have them be as balanced as possible, as they often have other commitments and responsibilities, which can be affected by having to participate in too many thesis defences.

\begin{figure}[H]
    \centering
    \includegraphics[scale=0.4]{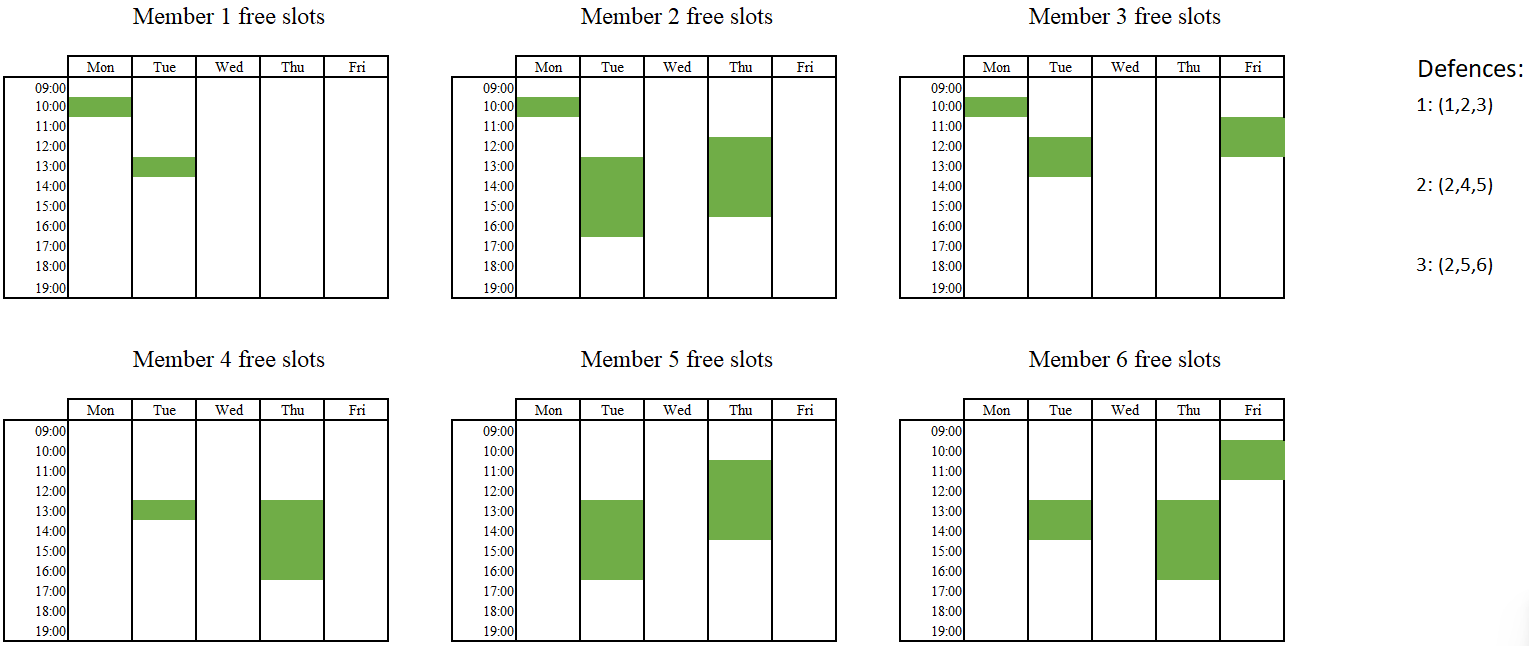}
    \caption{Time slot availability of six members}
    \label{avail}
\end{figure}

Assuming this assignment is already done, a time slot must be chosen for each defence. This is a time-consuming and demanding task for all of the people involved. When done by hand, all members' time availability must be gathered (Figure \ref{avail}), for example, using doodles or communicating \textit{via} e-mail. This is not always a simple task, as some members can forget or take a long time to respond and the person responsible must repeatedly contact some of them. Moreover, sometimes members answer \textit{via} different communication platforms, which ends up making this a somewhat chaotic process.

\begin{figure}[H]
    \centering
    \includegraphics[scale=0.4]{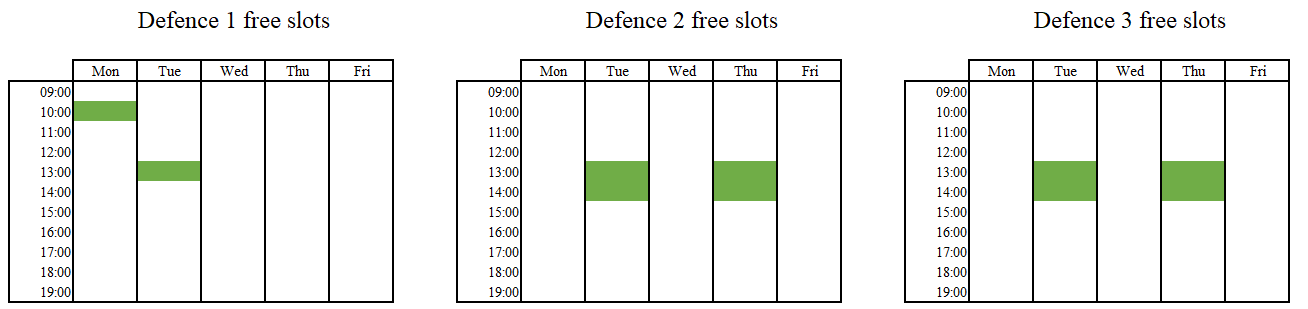}
    \caption{Feasible slots for each defence}
    \label{avail-d}
\end{figure}

Given these availabilities, it is necessary to determine the time slots that align for all committee members assigned to each defence, to ascertain when it is feasible to schedule them (Figure \ref{avail-d}). This becomes quite a daunting task after a certain number of defences. In some cases, there is no available time slot and the members must communicate and agree on some form of compromise. Other times there seems to be an available time slot for two defences, but the same member would have two assignments simultaneously, which is impossible.

Besides being an arduous task, scheduling thesis defences often leads to unsatisfactory timetables for the committee members. Their preferences over some time slots are not always respected, and there are unfair choices which cause some members to have defences scheduled and spread out over several days in an attempt to reduce the number of days of other members.

Automating thesis defence scheduling in universities can improve the quality of the defence assessment and make the process less burdensome for those tasked with its resolution and for the committee members whose daily tasks are impacted by this additional assignment. Studying this type of educational scheduling problem can enhance not just the operational effectiveness of universities, but also serve as a model for organisations seeking to address such inefficiencies and work towards streamlining and automating their procedures. This endeavour ensures not only improved performance in relevant KPIs but also liberates personnel to allocate their time to other tasks.

\section{The Monolithic Multi-Objective Mixed Integer Linear Programming Model}\label{s-mono}
\noindent This section presents the monolithic multi-objective mixed-integer linear programming model.

\subsection{Indices and sets}\label{SS-Model-Sets}
\noindent This subsection presents the indices and sets. It is important to highlight that although all entities start with a zero value, it is never employed to represent an actual object. Nevertheless, it is required to depict the nonexistence of these objects in certain constraints and objective functions. For instance, even in the absence of a committee member labelled as zero, we still require the capacity to quantify the presence of zero committees. 

\begin{enumerate}
    \item \textit{Indices.}
\begin{itemize}[label={--}]
    \item $i = 0,\ldots, n_i$, are the indices associated with master's thesis defence committee members;
    \item $j = 0,\ldots, n_j$, are the indices associated with the master's thesis defences;
    \item $t = 0,\ldots, n_t$, are the indices associated with the role of the committee members;
    \item $k = 0,\ldots, n_k$, are the indices associated with days;
    \item $\ell = 0,\ldots, n_\ell$, are the indices associated with the available hour slots in each day;
    \item $p = 0,\ldots, n_p$, are the indices associated with the available rooms;
    \item $q = 0,\ldots, n_q$, are the indices associated with research subjects.  
\end{itemize}
\item \textit{Sets.}

\begin{itemize}[label={--}]
            \item $A_{jt} = \{1,\ldots,i, \ldots,n_i\}$, is the set of committee members, $i$, that are eligible to be assigned to a role, $t$, in a defence, $j$, for all $j = 1,\ldots, n_j$, $t = 1,\ldots, n_t$;
             \item $A_{i} = \{(1,1),\ldots,(jt), \ldots,(n_j,n_t)\}$, is the set of ordered pairs of defence, $j$, and role, $t$, that a committee member, $i$, is eligible to be assigned to, for all $i = 1,\ldots, n_i$;

             \item $B_{i} = \{(1,1),\ldots,(k\ell), \ldots,(n_k,n_\ell)\}$, is set of available time slots, ($k\ell$), for each committee member, $i$, for all $i = 1,\ldots, n_i$;
            \item $C^{m}_{i} = \{1,\ldots,q, \ldots,n_q\}$, is the set of research subjects $q$, that a committee member, $i$, is an expert in, for all $i = 1,\ldots, n_i$;
             \item $C^{d}_{j}= \{1,\ldots,q, \ldots,n_q\}$, is the set of research subjects, $q$, that a defence, $j$, explores, for all $j = 1,\ldots, n_j$.
        \end{itemize}

\end{enumerate}
\subsection{Parameters}\label{SS-Model-Parameters}
\noindent The parameters for the model are presented in this subsection. 
        \begin{itemize}[label={--}]
            \item $a_{ik\ell} \in \mathbb{N}_0$, is the penalty that each committee member, $i$, holds for every time slot ($k\ell$), for all $i = 1,\ldots, n_i$, $k = 1,\ldots, n_k$, $\ell = 1,\ldots, n_\ell$;
            
            \item $d \in \mathbb{N}$, is the duration of a defence in hour slots.

        \end{itemize}

\subsection{Variables}\label{SS-Model-Variables}
\noindent This subsection presents the definition of the decision and auxiliary variables. 

\begin{enumerate}
    \item \textit{Decision variables.}
        \begin{itemize}[label={--}]
            \item $x_{ijtk\ell p}\in \{0,1\}$, is $1$ if an eligible committee member, $i$, is assigned to a thesis defence, $j$, performing a role, $t$, for which they are available, $(jt)\in A_i$. In an available time-slot ($k\ell$)$\in B_i$,  and in a room, $p$; and $0$ otherwise, for all $i = 1,\ldots, n_i$, $(jt) \in A_i$, $(k\ell)\in B_i$, and $p = 1,\ldots, n_p$.
        \end{itemize}
     \item \textit{Auxiliary variables.}
        \begin{itemize}[label={--}]  
            \item $y_{jk\ell p}\in \{0,1\}$, is 1 if the thesis defence $j$ is scheduled for day $k$, at hour $\ell$ in room $p$; and $0$ otherwise, for all $j = 1,\ldots, n_j$, $k = 1,\ldots, n_k$, $\ell = 1,\ldots, n_\ell$, and $p = 1,\ldots, n_p$; 
            \item$y_{ijk}\in \{0,1\} $, is $1$ if committee member $i$ is assigned to $j$ defences on day $k$; and $0$ otherwise, for all $i = 1,\ldots, n_i$, $j = 0, \ldots, n_j$ and $k = 1,\ldots, n_k$; 
            \item $w^{l}_{ij}\in \{0,1\} $, is $1$ if committee member $i$ is assigned to $j$ thesis defences; and $0$ otherwise, for all $i = 1,\ldots, n_i$, and $j = 0,\ldots, n_j$; 
            \item $ w^{d}_{ik}\in \{0,1\}$, is $1$ if committee member $i$ is assigned to thesis defences in $k$ days; and $0$ otherwise, for all $i = 1,\ldots, n_i$, and $k = 0,\ldots, n_k$.

        \end{itemize}
\end{enumerate}

\subsection{Objective functions}\label{SS-Model-Objectives}
\noindent This subsection presents the objective functions. Two points of view are considered. The first, committee assignment quality, includes the objectives that are determined by the selection of committee members for each defence. The second one, schedule quality, is related to the assignment of defences to certain time slots. However, additional constraints and variables must be defined before addressing their mathematical expressions. Subsection \ref{SS-Model-BackObjectives} properly defines the objective functions.

    \begin{enumerate}
        \item {\textit{Point of view of committee assignment quality.} This point of view includes the following objectives. }
            \begin{enumerate}
                \item \textit{Balance workloads}. The workload of a member is the number of committees they are assigned to. To ensure a balanced and fair assignment, the square of this number is to be minimised.

               \item\textit{Maximise committee member suitability}. Committee member suitability is measured as the number of research subjects each committee member has in common with their assigned defences. To ensure that the committees are composed of experts in the subjects of the defences, this number is to be maximised.
            \end{enumerate}
        \item \textit{Point of view of schedule quality.} This point of view includes the following objectives. 
            \begin{enumerate}

            \item \textit{Minimise the non-satisfaction of time slot preferences}. A penalty value is applied whenever a committee member is assigned an undesirable time slot. The penalty values are to be minimised.
            \item \textit{Minimise committee days}. Committee days are the number of days a committee member is scheduled to attend a defence. By minimising the square of this number, we also promote the fairness of the assignments.
          
        \end{enumerate}
    \end{enumerate}

\subsection{Constraints}\label{SS-Model-Constraints}
\noindent This subsection presents the constraints. We divide them into two categories. The first concerns the constraints that define the feasible region and the second the constraints that define the value for the auxiliary variables used in the objective functions

\begin{enumerate}
    \item \textit{Feasibility constraints} These constraints define the feasible region.
        \begin{enumerate}
            \item \textit{Complete committee definition}. A complete committee must include $n_t$ assignments for a defence, $j$, all with a different appointed role, $t$, in the same slot (day $k$, hour $\ell$, and room $p$). Let us note that an assignment for a defence can only be made if a member, $i$, is eligible to perform the respective role, \textit{i.e.}, $i\in A_{jt}$. Moreover, the defence can only be assigned to a time slot if its members are available for that slot, $(k\ell) \in B_{i}$.

                   \begin{equation}\label{Const - comm size - eq1}
                      \displaystyle \sum_{i \in A_{jt}} x_{ijtk\ell p} =  y_{jk\ell p},  \;\;\;  j = 1, \ldots, n_j,\; t = 1, \ldots, n_t, \; (k\ell) \in B_i,\; p = 1, \ldots, n_p.
                \end{equation}

            \item \textit{Single committee assignment}. {Each defence, $j$, can only be assigned to one committee and appointed one slot (day $k$, hour $\ell$, and room $p$). Thus, in this constraint, we state that for a defence, $j$, the number of complete committees assigned to it is $1$.}
  
                \begin{equation}\label{const - defences scheduled once - eq 4}
                    \displaystyle \sum_{k = 1}^{n_k}\sum_{\ell = 1}^{n_\ell}\sum_{p = 1}^{n_p} y_{jk\ell p} = 1, \;\;\;  j = 1, \ldots, n_j.
                \end{equation}
           
            \item \textit{Committee member assignment juxtaposition}. A committee member, $i$, cannot be assigned to more than one defence, $j$, at the same time. This constraint guarantees that a committee member is not assigned more than one role in a defence, as that would lead to multiple assignments in the same time slot. Let us note that assignments are only defined for pairs of defence and role ($jt$), if a member is eligible to perform the role, $t$, in that defence, $j$. Moreover, a member can only be assigned to a time slot if they are available for that slot, $(k\ell) \in B_{i}$.
                \begin{equation}\label{const - member juxtaposition}
                    \displaystyle \sum_{\{jt:(jt)\in A_{i}\}}\;\sum_{\ell = \overline{\ell},(k\ell) \in B_{i}}^{\overline{\ell} + d - 1}\;\sum_{p = 1}^{n_p} x_{ijtk\ell p} \leqslant 1, \;\;\; i = 1, \ldots, n_i,\; k = 1, \ldots, n_k,\; \overline{\ell} = 1, \ldots, n_\ell - d + 1.
                \end{equation}
            \item \textit{Room assignment juxtaposition}. A room, $p$, can only hold one defence, $j$, at a time ($k\ell$). 
                \begin{equation}\label{const - room capacity}
                    \displaystyle \sum_{j = 1}^{n_j}\sum_{\ell = \overline{\ell}}^{\overline{\ell} + d - 1} y_{jk\ell p} \leqslant 1, \;\;\; k = 1, \ldots, n_k,\; \overline{\ell} = 1, \ldots, n_\ell - d + 1,\; p = 1, \ldots, n_p.
                \end{equation}
        \end{enumerate}

        \item \textit{Auxiliary variable constraints.} These constraints define the values for the auxiliary variables necessary for the objective functions.
            \begin{enumerate}
                \item{ \textit{Workload definition.} The workload for a committee member, $i$, is the number, $j$, of defences they are assigned to. It is represented through a variable, $w^{l}_{ij}$, which takes the value 1 if a committee member, $i$, is assigned to a number, $j$, of defences, and 0 otherwise. }
                    \begin{equation}\label{const - wl def - eq1}
                        \displaystyle \sum_{j = 0}^{n_j} jw^{l}_{ij} =  \sum_{\{jt:(jt) \in A_i\}}\sum_{\{k\ell:(k\ell)\in B_i\}}\sum_{p = 1}^{n_p} x_{ijtk\ell p}, \;\;\; i = 1, \ldots, n_i.
                    \end{equation}
                    \begin{equation}\label{const - wl def - eq2}
                        \displaystyle \sum_{j = 0}^{n_j} w^{l}_{ij} = 1, \;\;\; i = 1, \ldots, n_i.
                    \end{equation}
                \item {\textit{Committee days definition.} A committee day occurs when a committee member, $i$, has a defence, $j$, scheduled. To represent these days, a variable,  $ w^{d}_{ik}$, takes the value $1$ if a committee member, $i$, has defences scheduled on a number, $k$, of days, and 0 otherwise. }

                \begin{equation}\label{const - comm days - eq 1}
                    \displaystyle\sum_{j = 0}^{n_j} jy_{ijk} =  \sum_{\{jt:(jt) \in A_{i}\}}\;\sum_{\{\ell:(k\ell) \in B_i\}}\;\sum_{p = 1}^{n_p} x_{ijtk\ell p}, \;\;\; i = 1, \ldots, n_i, \; k = 1, \ldots, n_k.
                \end{equation}
                \begin{equation}\label{const - comm days - eq 2}
                    \displaystyle \sum_{j = 0}^{n_j} y_{ijk} = 1, \;\;\; i = 1,\ldots, n_i, \; k = 1, \ldots, n_k.
                \end{equation}
                \begin{equation}\label{const - comm days - eq 4}
                    \displaystyle \sum_{k = 0}^{n_k} k w^{d}_{ik} = \sum_{j = 1}^{C^{m}_{i}}\sum_{k = 1}^{n_k} y_{ijk} , \;\;\; i = 1, \ldots, n_i.
                \end{equation}

                \begin{equation}\label{const - comm days - eq 5}
                    \displaystyle \sum_{k = 0}^{n_k}  w^{d}_{ik} = 1, \;\;\; i = 1, \ldots, n_i.
                \end{equation} 

            \end{enumerate}
    \end{enumerate}

\subsection{Back to the objective functions}\label{SS-Model-BackObjectives}
\noindent This subsection revisits and adequately defines the objective functions, now that we have presented all the necessary constraints.
    \begin{enumerate}
        \item {\textit{Point of view of committee assignment quality}. This point of view includes the following objectives.}
            \begin{enumerate}
                \item \textit{Minimise unfair workloads}. Minimising an exponential penalty promotes a fair workload distribution. The variable, $w^{l}_{ij}$, takes the value 1 if a committee member, $i$, is assigned to a number, $j$, of defences, and 0 otherwise. The linearity of the model is respected by multiplying $w^{l}_{ij}$ by $j^2$.
                    \begin{equation}\label{ob1}
                        {\displaystyle \min z_1(w) = \sum_{i = 1}^{n_i}\sum_{j = 1}^{n_j}j^{2}w^{l}_{ij}}.
                    \end{equation}
                \item \textit{Maximise committee member suitability}. Maximising the assignments of committee members, $i$, to defences, $j$, that they share research subjects with, $q \in C^{m}_{i} \cap C^{d}_{j}$, promotes the suitability of the committees to assess the defences.
                    \begin{equation}
                        {\displaystyle \max z_2(x) = \sum_{i = 1}^{n_i} \sum_{\{jt:(jt)\in A_{i}\}}\sum_{\{k,l:(k\ell) \in B_{i}\}}\sum_{p = 1}^{n_p}\sum_{\{q:q \in C^{m}_{i} \cap C^{d}_{j}\}}  x_{ijtk\ell p}}.
                    \end{equation}
        \end{enumerate}
        \item \textit{Point of view of schedule quality}. This point of view includes the following objectives.
            \begin{enumerate}

                \item \textit{Minimise the non-satisfaction of time slot preferences}. Minimising the assignments that are scheduled for time slots, ($k\ell$), that have a penalty associated, $a_{ik\ell}>0$, promotes the reduction of such assignments.
                    \begin{equation}
                        {\displaystyle \min z_3(x) = \sum_{i = 1}^{n_i} \sum_{\{jt:(jt)\in A_{i}\}}\sum_{\{k,l:(k\ell) \in B_{i}\}}\sum_{p = 1}^{n_p} a_{ik\ell} x_{ijtk\ell p}}.
                    \end{equation}
                \item \textit{Minimise committee days}. Minimising an exponential penalty promotes the minimisation and fair distribution of the committee days. The variable, $ w^{d}_{ik}$, takes the value 1 if a committee member, $i$, is assigned to defences in a number, $k$, of days, and 0 otherwise. The linearity of the model is respected by multiplying $w_{ik}$ by $k^2$.
                    \begin{equation}
                        {\displaystyle \min z_4({w}^d) = \sum_{i = 1}^{n_i}\sum_{k = 1}^{n_k}k^{2} w^{d}_{ik}}.
                    \end{equation}

            \end{enumerate}
    \end{enumerate}

\section{The Decomposed Multi-Objective Mixed Integer Linear Programming Model}\label{s-dec}
\noindent This section presents the decomposed multi-objective mixed-integer linear programming model. It includes two submodels to be solved sequentially. The first submodel assigns committees to defences and the second schedules those defences for a certain time slot, considering the committee configurations defined in the first submodel.

The indices and parameters introduced in the previous section are also the ones considered for the decomposed model. The second submodel is similar to the monolithic model, with two main differences. Since the committees are already defined, no committee assignment quality objective is assessed. Moreover, the set of committee members, $i$, that are eligible to perform a certain role, $t$, in a defence, $j$, $A_{jt}$, is comprised of a single member for each defence and role, as the committees are already defined. Consequently, each member is only eligible for the pairs of defences and roles, $A_i$, they are pre-assigned to. As these are the only distinctions between the second submodel and the monolithic model, the remainder of this section is dedicated to the first submodel, which defines the committee configurations.

\subsection{Variables}

\noindent This subsection presents the definition of the decision and auxiliary variables for the first submodel.
\begin{enumerate}

\item \textit{Decision variables.}
        \begin{itemize}[label={--}]
            \item $ x_{ijt}\in \{0,1\}$, is $1$ if an eligible committee member, $i$, is assigned to a thesis defence, $j$, performing a role, $t$, for which they are available, $(jt)\in A_i$; and $0$ otherwise, for all $i = 1,\ldots, n_i$, $(jt) \in A_i$.
        \end{itemize}
     \item \textit{Auxiliary variables.}
        \begin{itemize}[label={--}]  
            \item $w^{l}_{ij}\in \{0,1\} $, is $1$ if committee member $i$ is assigned to $j$ thesis defences; and $0$ otherwise, for all $i = 1,\ldots, n_i$, and $j = 0,\ldots, n_j$; 
            \item $w_{jk\ell} \in \{0,\ldots,n_t\}$, is the number of committee members assigned to a defence, $j$, that is available for a time slot, ($k\ell$), for all $j= 1,\ldots,n_j$, $k=1,\ldots,n_k$, $\ell= 1,\ldots,n_\ell$;
            \item $ w_{jtk\ell} \in \{0,1\}$, is 1 if the number of committee members assigned to a defence, $j$, that is available for a time slot, ($k\ell$), is equal to $t$; and $0$ otherwise, for all $j= 1,\ldots,n_j$, $t = 0,\ldots,n_t$, $k=1,\ldots,n_k$, $\ell= 1,\ldots,n_\ell$.
        \end{itemize}
\end{enumerate}

\subsection{Objective functions}

\noindent The first submodel aims to define the committees. Accordingly, the objectives that model the committee assignment quality are included. In this stage, the defences are not assigned to any time slot. Nonetheless, we want to choose committees which have the potential to have good performances regarding this point of view. Thus, we have defined a proxy objective that represents this potential.

  \begin{enumerate}
        \item {\textit{Point of view of committee assignment quality.} This point of view includes the following objectives. }
            \begin{enumerate}
                \item \textit{Minimise unfair workloads}. Similar to the monolithic model.

               \item\textit{Maximise committee member suitability}. Similar to the monolithic model.
            \end{enumerate}
        \item \textit{Point of view of schedule quality.} This point of view includes the following objectives. 
            \begin{enumerate}

            \item \textit{Maximise number of available time slots}. An available time slot for a defence is defined as a time slot where all its assigned members are available. The total number of available time slots provides a proxy for the performance of the schedule quality objectives.          
           
        \end{enumerate}
\end{enumerate}
\subsection{Constraints}

\noindent This subsection presents the constraints for the first submodel. We divide them into two categories. The first concerns
the constraints that define the feasible region and the second the constraints that define the value for the
auxiliary variables used in the objective functions.

\begin{enumerate}

 \item \textit{Feasibility constraints} These constraints are used to define the feasible region.
        \begin{enumerate}
            \item \textit{Complete committee definition}. A complete committee must include $n_t$ assignments for a defence, $j$, all with a different appointed role, $t$. Let us note that an assignment for a defence can only be made if a member, $i$, is eligible to perform the respective role, \textit{i.e.}, $i\in A_{jt}$.

                   \begin{equation}\label{Const - comm size - eq1}
                      \displaystyle \sum_{i \in A_{jt}}  x_{ijt} =  1,  \;\;\;  j = 1, \ldots, n_j,\; t = 1, \ldots, n_t. 
                \end{equation}

            \item \textit{Single committee member assignment}. Each committee member, $i$, can only be assigned once to a defence, $j$, they are eligible to be assigned to, $(jt)\in A_i$.
                \begin{equation}\label{Const - comm size - eq1}
                      \displaystyle \sum_{\{jt:(jt)\in A_i\}}  x_{ijt} \leqslant  1,  \;\;\;  i = 1, \ldots, n_i,\;  j = 1, \ldots, n_j.
                \end{equation}
            \item \textit{Time slot availability}. Each defence, $j$, has at least one time slot it can be assigned to. The variable, $ w_{jtk\ell}$, is 1 if the number of committee members assigned to a defence, $j$, that is available for a time slot, ($k\ell$), is equal to $t$. Accordingly, for each defence, $j$, there must be at least one time slot, ($k\ell$), such that $n_t$ of its committee members are available.
            \begin{equation}
                \displaystyle \sum_{k=1}^{n_k}\sum_{\ell=1}^{n_\ell} w_{jn_tk\ell}\geqslant 1, \;\;\; j = 1,\ldots,n_j.
            \end{equation}
           
           \end{enumerate}     

        \item \textit{Auxiliary variable constraints.} These constraints define the values for the auxiliary variables necessary for some of the objective functions and feasibility constraints.
            \begin{enumerate}
                \item{ \textit{Workload definition.} The workload for a committee member, $i$, is the number, $j$, of defences they are assigned to. It is represented through a variable, $w^{l}_{ij}$, which takes the value 1 if a committee member, $i$, is assigned to a number, $j$, of defences, and 0 otherwise. }
                    \begin{equation}\label{const - wl def - eq1}
                        \displaystyle \sum_{j = 0}^{n_j} jw^{l}_{ij} =  \sum_{\{jt:(jt) \in A_i\}}  x_{ijt}, \;\;\; i = 1, \ldots, n_i.
                    \end{equation}
                    \begin{equation}\label{const - wl def - eq2}
                        \displaystyle \sum_{j = 0}^{n_j} w^{l}_{ij} = 1, \;\;\; i = 1, \ldots, n_i.
                    \end{equation}
                \item \textit{Time slot availability definition.} Variable $w_{jk\ell}$ represents the number of committee members assigned to a defence, $j$, that is available for a time slot, ($k\ell$). Constraint \ref{tsd1} defines this variable. For each defence, $j$, and time slot, ($k\ell$), the variable is defined as the sum of the assignments that have a committee member, $i$, available for that slot, $(k\ell)\in B_i$. Variable $ w_{jtk\ell}$, is 1 if the number of committee members assigned to a defence, $j$, that is available for a time slot, ($k\ell$), is  $t$. Constraints \ref{tsd2} and \ref{tsd3} are used to define this variable.

                \begin{equation}\label{tsd1}
                    \displaystyle \sum_{t = 1}^{n_t}\sum_{\{i:i\in A_{jt}, (k\ell)\in B_i\}} x_{ijt}=w_{jk\ell},\;\;\; j = 1,\ldots,n_j,\;k = 1,\ldots,n_k,\;\ell=1,\ldots,n_\ell.
                \end{equation}

                 \begin{equation}\label{tsd2}
                    \displaystyle \sum_{t = 0}^{n_t}t w_{jtk\ell}=w_{jk\ell},\;\;\; j = 1,\ldots,n_j,\;k = 1,\ldots,n_k,\;\ell=1,\ldots,n_\ell.
                \end{equation}

                 \begin{equation}\label{tsd3}
                    \displaystyle \sum_{t = 0}^{n_t} w_{jtk\ell}=1,\;\;\; j = 1,\ldots,n_j,\;k = 1,\ldots,n_k,\;\ell=1,\ldots,n_\ell.
                \end{equation}
            \end{enumerate}
    \end{enumerate}

\subsection{Back to the objective functions}\label{S-Model-BackObjectives}
\noindent This subsection revisits and adequately defines the objective functions, now that we have presented all the necessary constraints.
    \begin{enumerate}
        \item {\textit{Point of view of committee assignment quality}. This point of view includes the following objectives.}
            \begin{enumerate}
                \item \textit{Minimise unfair workloads}. Minimising an exponential penalty promotes a fair workload distribution. The variable, $w^{l}_{ij}$, takes the value 1 if a committee member, $i$, is assigned to a number, $j$, of defences, and 0 otherwise. The linearity of the model is respected by multiplying $w^{l}_{ij}$ by $j^2$.
                    \begin{equation}\label{ob1}
                        {\displaystyle \min z_1(w) = \sum_{i = 1}^{n_i}\sum_{j = 1}^{n_j}j^{2}w^{l}_{ij}}.
                    \end{equation}
                \item \textit{Maximise committee member suitability}. Maximising the assignments of committee members, $i$, to defences, $j$, that they share research subjects with, $q \in C^{m}_{i} \cap C^{d}_{j}$, promotes the suitability of the committees to assess the defences.
                    \begin{equation}
                        {\displaystyle \max z_2(x) = \sum_{i = 1}^{n_i} \sum_{\{jt:(jt)\in A_{i}\}}\sum_{\{q:q \in C^{m}_{i} \cap C^{d}_{j}\}}   x_{ijt}}.
                    \end{equation}
        \end{enumerate}
        \item \textit{Point of view of schedule quality}. This point of view includes the following objective.
            \begin{enumerate}

                \item \textit{Maximise the available time slots}. An available time slot, ($k\ell$), for a defence, $j$, is defined as one where all three of the committee members assigned to that defence are available, $ w_{jn_tk\ell}=1$. Maximising the available time slots promotes the flexibility of the assignments in the next stage and provides a proxy for the performance of the objectives from the point of view of schedule quality.
                    \begin{equation}
                        {\displaystyle \max z_5(x) = \sum_{j = 1}^{n_j} \sum_{k=1}^{n_k}\sum_{\ell=1}^{n_\ell} w_{jn_tk\ell}}.
                    \end{equation}

            \end{enumerate}
    \end{enumerate}

\section{Augmented $\epsilon$-constraint methods}\label{s-augm}

\noindent Augmented $\epsilon$-constraint methods are used to solve both the monolithic problem and the second submodel of the decomposed problem. This section provides a theoretical background regarding these methods. The practical aspects are further discussed in Section \ref{s-AlgF}.

\subsection{Concepts, definitions, and notation}

\noindent Consider the Multi-Objective Mixed Integer Linear Programming (MOMILP):

\begin{equation}
    \begin{array}{l}
        \max \;z_1(x),   \\
        \;\;\;\;\; \vdots \\
        \max \;z_i(x), \\
        \;\;\;\;\;  \vdots\\
        \max\; z_{n_z}(x),\\\\
        \text{subject to:}\\
        x \in X
    \end{array}
\end{equation}

Where $x=(x_1,\ldots,x_j,\ldots,x_{n_x})$ is The vector of integer and/or continuous decision variables, the feasible region in the decision space is $X$, and the vector of linear objective functions is $z_i(x),\; i=1,\ldots,n_z$. The image of $X$ considering the objective functions defines the feasible region, $Z$, in the objective space.

A solution or outcome vector $z^\prime$ in the objective space dominates another solution $z^{\prime\prime}\neq z^{\prime}$ if and only if $z^{\prime}_i \geq z^{\prime\prime}_i$ for all $i=1,\ldots,n_z$, with at least one strict inequality, i.e., $z^{\prime}_i > z^{\prime\prime}_i$ for some $i$.
A feasible solution, $\bar{z} \in Z$, is \textit{non-dominated} if and only if there is no other feasible solution, $z \in Z$, such that $z$ dominates $\bar{z}$. The set of all non-dominated solutions is known as the \textit{non-dominated front}, $N$, or \textit{Pareto} front. The inverse image of a non-dominated solution, $\bar{x} = f^{-1}(\bar{z})$, is called an \textit{efficient solution} in the decision space.

Our objective in this paper is to identify a subset of the non-dominated front, $S\subseteq N$. For this purpose, a well-known scalarisation technique based on solving a sequence of $\epsilon-$constraint problems is used. These problems can be stated as follows:

\begin{equation}
    \begin{array}{l}
        \max \;z_1(x),  \\
       \\
        \text{subject to:}\\
        x \in X,\\
        z_i(x) \geqslant \epsilon_i, \; i = 2, \ldots, n_z.
    \end{array}
\end{equation}

A single objective function, arbitrarily defined as $z_1(x)$, is maximised, while the others are included in constraints, setting lower bounds, $\epsilon_i$, for objective functions $z_i(x)$, with $i=2,\ldots,n_z$. Different non-dominated solutions are obtained by setting different values for the lower bounds, $\epsilon_i$.

\subsection{The augmented $\epsilon$-constraint by \cite{MavrotasFlorios2013}}

\noindent The advantages the augmented $\epsilon$-constraint method \citep{MavrotasFlorios2013} has in comparison with other multi-objective optimisation algorithms are that all of the solutions found are guaranteed to be non-dominated and some lower or upper bound vectors, $\epsilon$, can be proven to be infeasible or leading to a previously found solution and can be skipped.
This is achieved by augmenting the integer objective function, $z_1(x)$, with a factor that includes the values of other objective functions, $z_i(x),\; i=2,\ldots,n_z$, guaranteeing that all the solution vectors found are non-dominated. Moreover, this factor must be small enough in such a way that the value of the integer objective function, $z_1(x)$, is not affected by it. An example of such a function can be stated as follows: 

\begin{equation}\label{aug_function}
    \displaystyle \max z_e(x) = z_1(x) + {n_z}^{-1}  \sum_{i=2}^{n_z} \frac{z_i(x) - z_i^{min}}{z_i^{max}-z_i^{min}}.
\end{equation}

This function includes the integer objective function, $z_1(x)$, and the augmentation factor. This factor is the sum of the \textit{ratios} of the difference between the value of the remaining objective functions, $z_i(x),\; i=2,\ldots,n_z$, and the minimum value of each objective function, $z_i^{min},\; i=2,\ldots,n_z$, and the difference between the maximum value of the objective functions, $z_i^{max},\; i=2,\ldots,n_z$, and their minimum values, $z_i^{min},\; i=2,\ldots,n_z$. Let us note that the maximum value each of these \textit{ratios} can take is 1. Accordingly, the maximum value the sum can take is $n_z-1$. Thus, by dividing their sum by $n_z$, it is guaranteed that the value of the factor is never large enough to impact the value of the integer objective function, $z_1(x)$. Different procedures can be used to compute the values for the vectors, $z^{min}$ and $z^{max}$. Our procedure is presented in Section \ref{initi}.

The pseudo-code for an augmented $\epsilon$-constraint method which uses an objective function like Expression \eqref{aug_function} is provided in Algorithm \ref{alg_econst}. The inputs are the vector of minimum values for each of the bounded objectives, $z^{min} =(z_2^{min},\ldots,z_i^{min},\ldots,z_{n_z}^{min})$, the vector of maximum values, $z^{max} =(z_2^{max},\ldots,z_i^{max},\ldots,z_{n_z}^{max})$, the vector of lower bound increments, $\epsilon' =(\epsilon' _2,\ldots,\epsilon' _i,\ldots,\epsilon' _{n_z})$, and the problem to be optimised, $P$. The output of this algorithm is a set of non-dominated solutions, $S$.

In the lower bound vector, $\epsilon$, the starting point is the minimum values vector, $z^{min}$. And in the last iteration, it is the maximum values vector, $z^{max}$. In each iteration, the value of one of the lower bounds, $\epsilon_i$, is incremented by a factor, $\epsilon'_i$. If all of the objective functions are integers and this increment factor is 1, all of the non-dominated solutions between the minimum and maximum vectors can be found.

The problem is optimised in each iteration considering the lower bounds vector, $\epsilon$. However, some iterations can be skipped. There are two justifications to skip an iteration. The first is that there is a previously found solution, $s\in S$, where all of the bounded objectives have better or equal performances than the current lower bounds vector. This iteration is redundant as an equivalent solution will be found. The other justification is that a previous set of lower bounds, $\epsilon\in I$, was proven infeasible, and all were smaller or equal to the current set. This iteration is redundant as this set of lower bounds must be infeasible. Besides saving the results of feasible iterations, the sets of infeasible lower bounds must be kept in memory to do this skipping procedure.

\begin{algorithm}[h!]
\caption{Augmented $\epsilon$-constraint}\label{alg_econst}
	\begin{algorithmic}[1]\normalsize
	\State{\textbf{input:} $z^{min}$ (vector of minimum values), $z^{max}$ (vector of maximum values), $\epsilon'$ (vector of lower bound increments), $P$ (problem)};
	\State{\textbf{output:} $S$ (non dominated solutions)};
        \State{$S\leftarrow\{\}$;}
        \State{$I\leftarrow\{\}$;}
       \For{($\epsilon_2 = z_2^{min}; \epsilon_2 \leqslant z_2^{max}; +\epsilon'_2$)}
       \State{$\ldots$}
        \For{($\epsilon_i = z_i^{min}; \epsilon_i \leqslant z_i^{max}; +\epsilon'_i$)}
        \State{$\ldots$}
        \For{($\epsilon_{n_z} = z_{n_z}^{min}; \epsilon_{n_z} \leqslant z_{n_z}^{max}; +\epsilon'_{n_z}$)}

       \If{($not\_skip(\epsilon,S,I)$)}
  
       \State{$s\leftarrow optimise(P,\epsilon)$;}
       \If{($feasible(s)$)}
       \State{$S \leftarrow S \cup \{s\}$;}
       \Else
       \State{$I \leftarrow I\cup \{\epsilon\}$;}
       \EndIf
       \EndIf
        \EndFor
        \EndFor
        \EndFor
      
        \State{{\textbf{return}($\mbox{$S$}$);}}
        \end{algorithmic} 
	   
\end{algorithm}

\section{Genetic algorithms}\label{s-gen}

\noindent Genetic algorithms are population-based heuristic search methods. They are used to solve the first submodel of the decomposed problem, which aims to find committee compositions. This section presents the chromosome representation, the population initialisation procedure, the crossover and mutation operators, the crowding distance tournaments, and finally it addresses the genetic algorithms implemented for this work, NSGA-II and NSGA-III.

\subsection{Chromosome representation}\label{chromo}

\noindent In genetic algorithms, each solution, $s$, is called an individual, and its objective function value is used to determine its fitness. Each individual is defined by a chromosome from which its fitness can be obtained. The set of chromosomes considered in each iteration is called a population. 

\begin{figure}[h]
\centering
\scalebox{0.7}{

\begin{tikzpicture}
 \node[blockns] (1) {Role 1\\Member 2};
 \node[blockns, right = of 1, xshift=-1cm] (2) {Role 2\\Member 1};
 \node[blockns, right = of 2, xshift=-1cm] (3){Role 3\\Member 4};
 \node[block3, above = of 2, yshift=-0.5cm] (10) {Defence 1};
 \node[block0, right = of 3, xshift=-1cm] (4) {Role 1\\Member 1};
 \node[block0, right = of 4, xshift=-1cm] (5) {Role 2\\Member 2};
 \node[block0, right = of 5, xshift=-1cm] (6) {Role 3\\Member 5};
  \node[block3, above = of 5, yshift=-0.5cm] (11) {Defence 2};
 \node[blocks, right = of 6, xshift=-1cm] (7) {Role 1\\Member 2};
 \node[blocks, right = of 7, xshift=-1cm] (8) {Role 2\\Member 4};
  \node[blocks, right = of 8, xshift=-1cm] (9) {Role 3\\Member 3};
   \node[block3, above = of 8, yshift=-0.5cm] (12) {Defence 3};

   \draw (-1.1,-1.5)--(-1.1,1.5);
   \draw (5.5,-1.5)--(5.48,1.5);
   \draw (12.07,-1.5)--(12.07,1.5);
   \draw (18.67,-1.5)--(18.67,1.5);
\end{tikzpicture}}
   
    \caption{Chromosome representation for a problem with three defences, three roles, and five committee members.}
    \label{Chromo diagram}
\end{figure}
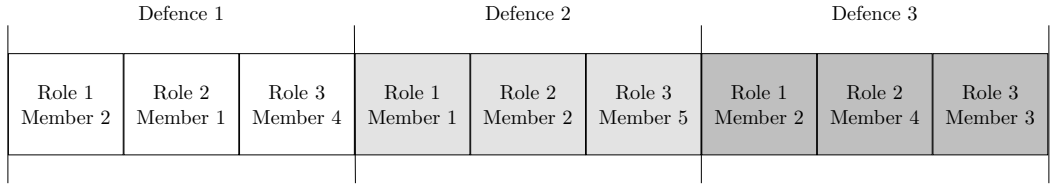

In our chromosome representation, for each defence, there are $n_t$ genes, one for each role in a committee. Thus, the chromosomes are of a fixed size, $n_jn_t$. Each gene represents a committee member, $i$, being assigned to a defence, $j$, and to a role, $t$. This representation means that the complete committee constraint is always fulfilled. Figure \ref{Chromo diagram} presents a graphical representation of a chromosome for a problem with three defences, three roles, and five committee members.

\subsection{Initialisation}\label{ini}

\noindent In genetic algorithms, the initialisation phase aims to generate new individuals, $s$, until the initial population, $S^{init}$, reaches a desired size, $n_{s}$. Algorithm \ref{alg_i} provides the pseudo-code for the initialisation of our genetic algorithms. In each iteration, a new individual, $s$, is added to the initial population, $S^{init}$. To generate this new individual, for every defence, $j$, a feasible committee is generated based on its eligible members, $A_{jt}$, and the available time slots for each member, $B_i$. This committee, $generate\_feasible\_committee(j,A_{jt},B_i)$, is then added, $\frown$, to the chromosome of the individual, $s$. The $length()$ function is used to check when the initial population reaches the desired size. It receives a set as input. Its output is the number of elements in that set. Here, this set is the initial population, $S^{init}$.

\begin{algorithm}[h!]
\caption{Initialise population}\label{alg_i}
	\begin{algorithmic}[1]\normalsize
	\State{\textbf{input:} $n_{s}$ (population size), $S^{init}$ (initial population), $A_{jt}$ (eligibility parameters), $B_{i}$ (availability parameters)};
	\State{\textbf{output:} $S^{init}$ (initial population)};
        
       \While{($length(S^{init})<n_{s}$)}

        \State{$s\leftarrow\{\}$;}
        \For{($j=1;j\leqslant n_j; +1$)}
         \State{$s\leftarrow s \frown generate\_feasible\_committee(j,A_{jt},B_i)$;}
        \EndFor
         \State{$S^{init}\leftarrow S^{init}\cup \{s\}$;}
        \EndWhile
      
        \State{{\textbf{return}($\mbox{$S^{init}$}$);}}
        \end{algorithmic} 
	   
\end{algorithm}

\begin{figure}[h!]
    \centering
    \includegraphics[scale=0.5]{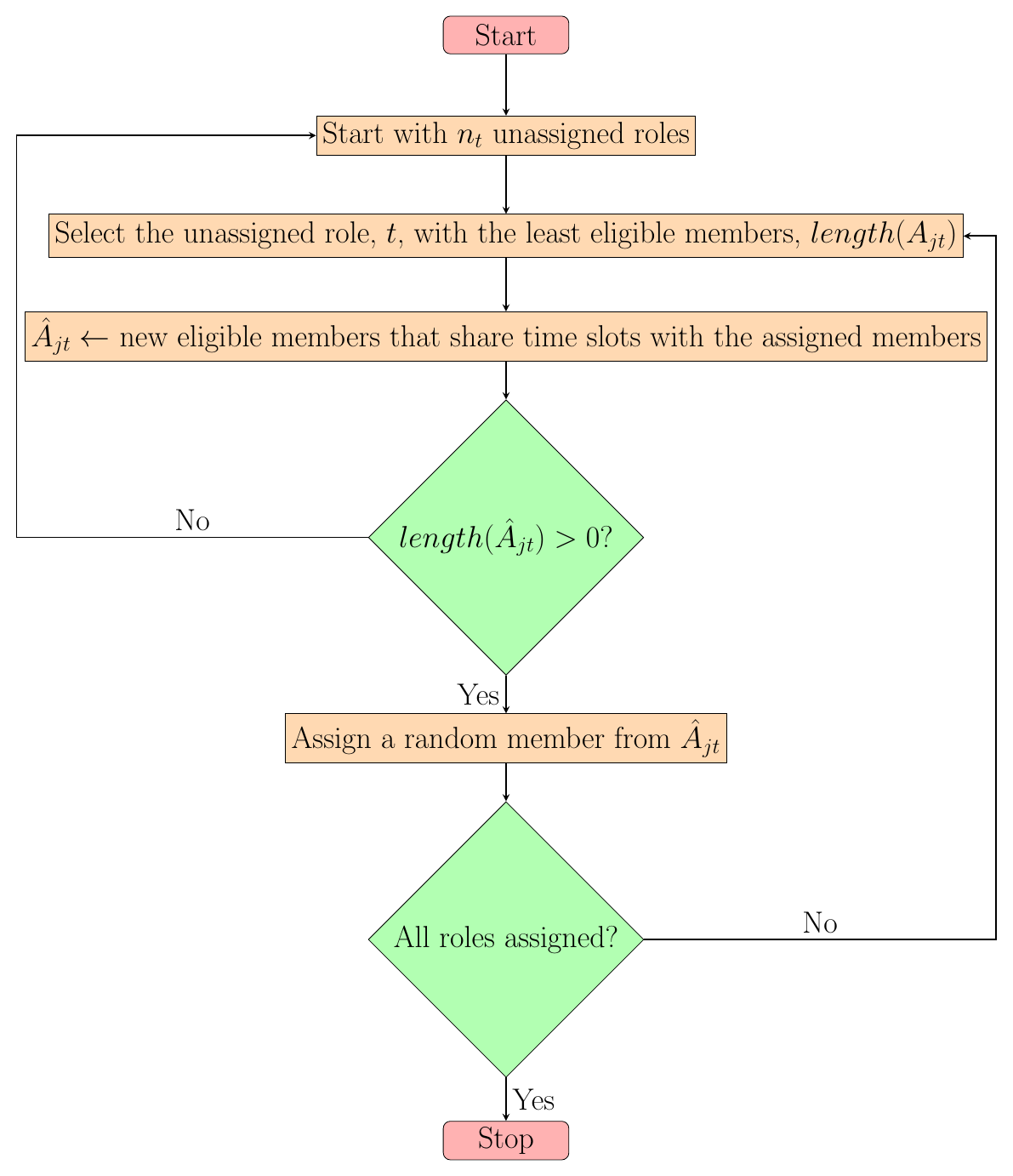}
    \caption{Generate feasible committee diagram}
    \label{genfeas}
\end{figure}

The pseudo-code diagram for generating each feasible committee is provided in Figure \ref{genfeas}. 

\subsection{Crossover}\label{cros}

\noindent Crossover is a genetic algorithm operator used to generate new individuals. Two parent chromosomes undergo crossover, creating a new chromosome based on a random recombination of their genes. The pseudocode for our crossover operator is given in Algorithm \ref{alg_cross}. The inputs are two parent chromosomes, $s^{p_1}$ and $s^{p_2}$, and the output is an offspring chromosome, $s^{off}$. For every defence, $j$, a random parent is chosen, $s^p$, and the $n_t$ genes of the parent representing the committee for the defence are passed down to the offspring chromosome, $s^{off} \frown s^p_j$. Entire committees are passed down to the offspring. This guarantees the feasibility of the offspring contingent on the parents being feasible. Contrarily, if we did the recombination based on single roles, this would not be guaranteed, as combinations with no available time slots or with duplicated members could be created. The $rand()$ function is used to choose which parent is randomly selected. It receives as input elements to be selected. Its output is the selected element. A probability for each element can be provided. If it is omitted it is uniform. Here, the input elements are the parent chromosomes, $s^{p_1}$ and $s^{p_2}$, and the probability is uniform.

\begin{algorithm}[h!]
\caption{Crossover}\label{alg_cross}
	\begin{algorithmic}[1]\normalsize
	\State{\textbf{input:} $s^{p_1}$ (parent 1), $s^{p_2}$ (parent 2)};
	\State{\textbf{output:} $s^{off}$ (offspring)};

        \State{$s^{off}\leftarrow\{\}$;}
        \For{($j=1;j\leqslant n_j; +1$)}
        \State{$s^p\leftarrow rand(s^{p_1},s^{p_2})$;}
         \State{$s^{off}\leftarrow s^{off} \frown s^p_j$;}
        \EndFor

        \State{{\textbf{return}($\mbox{$s^{off}$}$);}}
        \end{algorithmic} 
	   
\end{algorithm}

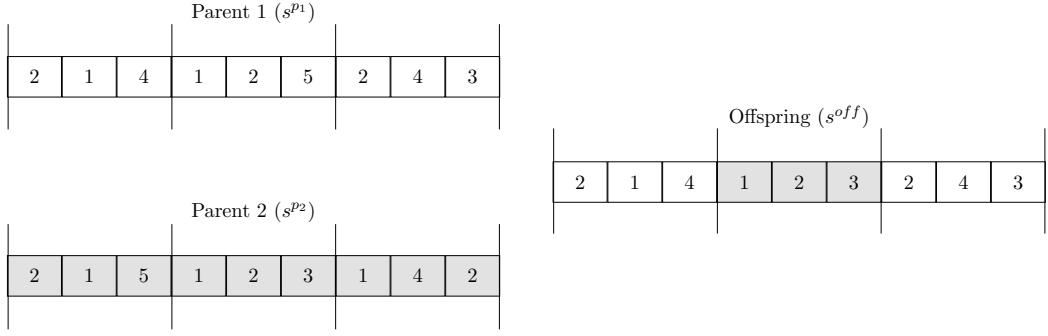
\begin{figure}[h]
\centering
\scalebox{0.7}{
\begin{tikzpicture}
 \node[blockns1] (1) {2};
 \node[blockns1, right = of 1, xshift=-1cm] (2) {1};
 \node[blockns1, right = of 2, xshift=-1cm] (3){4};

 \node[blockns1, right = of 3, xshift=-1cm] (4) {1};
 \node[blockns1, right = of 4, xshift=-1cm] (5) {2};
 \node[blockns1, right = of 5, xshift=-1cm] (6) {5};
  \node[block3, above = of 5, yshift=-0.5cm] (10){Parent 1 ($s^{p_1}$)};
 \node[blockns1, right = of 6, xshift=-1cm] (7) {2};
 \node[blockns1, right = of 7, xshift=-1cm] (8) {4};
  \node[blockns1, right = of 8, xshift=-1cm] (9) {3};

 \node[blocks1, below = of 1, yshift=-2cm] (11) {2};
 \node[blocks1, right = of 11, xshift=-1cm] (12) {1};
 \node[blocks1, right = of 12, xshift=-1cm] (13){5};

 \node[blocks1, right = of 13, xshift=-1cm] (14) {1};
 \node[blocks1, right = of 14, xshift=-1cm] (15) {2};
 \node[blocks1, right = of 15, xshift=-1cm] (16) {3};
  \node[block3, above = of 15, yshift=-0.5cm] (20){Parent 2 ($s^{p_2}$)};
 \node[blocks1, right = of 16, xshift=-1cm] (17) {1};
 \node[blocks1, right = of 17, xshift=-1cm] (18) {4};
  \node[blocks1, right = of 18, xshift=-1cm] (19) {2};

   \node[blockns1, right = of 9, yshift=-2cm] (21) {2};
 \node[blockns1, right = of 21, xshift=-1cm] (22) {1};
 \node[blockns1, right = of 22, xshift=-1cm] (23){4};

 \node[blocks1, right = of 23, xshift=-1cm] (24) {1};
 \node[blocks1, right = of 24, xshift=-1cm] (25) {2};
 \node[blocks1, right = of 25, xshift=-1cm] (26) {3};
  \node[block3, above = of 25, yshift=-0.5cm] (30){Offspring ($s^{off}$)};
 \node[blockns1, right = of 26, xshift=-1cm] (27) {2};
 \node[blockns1, right = of 27, xshift=-1cm] (28) {4};
  \node[blockns1, right = of 28, xshift=-1cm] (29) {3};

  \draw (-0.5,-1)--(-0.5,1);
  \draw (2.61,-1)--(2.61,1);
  \draw (5.72,-1)--(5.72,1);
  \draw (8.83,-1)--(8.83,1);

    \draw (-0.5,-4.8)--(-0.5,-2.8);
  \draw (2.61,-4.8)--(2.61,-2.8);
  \draw (5.72,-4.8)--(5.72,-2.8);
  \draw (8.83,-4.8)--(8.83,-2.8);

   \draw (9.86,-2.98)--(9.86,-0.98);
  \draw (12.97,-2.98)--(12.97,-0.98);
  \draw (16.08,-2.98)--(16.08,-0.98);
  \draw (19.19,-2.98)--(19.19,-0.98);
\end{tikzpicture}}
   
    \caption{Crossover representative diagram for a problem with three defences, three roles, and five committee members.}
    \label{Cross diagram}
\end{figure}

A representative diagram for a problem with three defences, three roles, and five committee members is provided in Figure \ref{Cross diagram}. For the first defence, the first parent is randomly selected to pass its genes, (2, 1, 4). For the second defence, the second parent is randomly selected and passes the genes, (1, 2, 3). For the third defence, the first parent is randomly selected and passes the genes, (2, 4, 3).

\subsection{Mutation}\label{mutation}

\noindent Mutation is a genetic algorithm operator that changes existing or newly generated individuals. If an individual undergoes mutation some of its genes are modified. The pseudocode for our mutation operator is given in Algorithm \ref{alg_mut}. The inputs are the initial individual, $s^{init}$, and the mutation probability, $m$, and the output is the mutated individual, $s^{mut}$. With a probability, $m$, an individual undergoes mutation. In that case, a new feasible committee is generated for a random defence, $j$, following the same protocol as previously depicted in Figure \ref{genfeas}.

\begin{algorithm}[h!]
\caption{Mutation}\label{alg_mut}
	\begin{algorithmic}[1]\normalsize
	\State{\textbf{input:} $s^{init}$ (initial individual), $m$ (mutation probability)};
	\State{\textbf{output:} $s^{mut}$ (mutated individual)};
         \State{$s^{mut}\leftarrow s^{init}$;}
        \If{($rand(1,\ldots,100)\leqslant m$ )}
        \State{$j\leftarrow rand(1,\ldots,n_j)$;}
        \State{$s^{mut}_{j}\leftarrow generate\_feasible\_committee(j, A_{jt}, B_i)$;}
        \EndIf

        \State{{\textbf{return}($\mbox{$s^{mut}$}$);}}
        \end{algorithmic} 
	   
\end{algorithm}

\subsection{Crowding distance tournament}\label{cdt}

\noindent The crowding distance measures how distant one point is to other points in the solution space. In our implementation of NSGA-II and NSGA-III, a tournament based on crowding distances is used to select which parent individuals are included in a mating pool of potential candidates for crossover. 

The crowding distance of a solution with respect to a group of solutions, $c(s,S)$, is defined in Equation \eqref{eq cd}. The closest point to a solution considering an objective function, $z_i$, with a greater or equal performance in that objective is denoted by $s_{i}^+$. The closest point considering an objective function, $z_i$, with a smaller or equal performance in that objective is denoted by $s_{i}^-$. For every objective function, $z_i$, the crowding distance is incremented by the difference between these two points. If there is no point with greater or equal or lower or equal performances reference values are used instead.

\begin{equation} \label{eq cd}
    \displaystyle c(s,S)=\sum_{i=1}^{n_z} \frac{z_i(s^+_i)-z_i(s^-_i)}{\max z_i(S) - \min z_i(S) } .
\end{equation}

\begin{figure}[H]
    \centering
\begin{tikzpicture}
	\begin{axis}[%
 xlabel={$z_1(x)$},
	ylabel={$z_2(x)$}
	={(xy): \thisrow{label}},%
	scatter/classes={%
		c={mark=triangle*,draw=black},
        d={mark=Mercedes star flipped,black},
        e={mark=asterisk,black}}    ]
	\addplot[scatter,only marks,%
		scatter src=explicit symbolic]%
	table[meta=label] {
x     y      label

384	203     c
388	196     c
390	189     c
400	182     c
416	175     c
436	147     c
502	143     c

368	197     d
370	182     d
372	179     d
376	173     d
378	161     d
380	158     d
382	155     d
386	149     d
390	146     d
396	141     d
402	137     d
412	132     d
420	125      d
488	121      d

338	145     e
340	139     e
342	135     e
344	129     e
346	123     e
348	120     e
350	118      e
354	112      e
358	108      e
362	104      e
374	98      e
378	94      e
392	90      e

	};
 \legend{ \tiny Non-dominated front 2 ($N_2$),\tiny  Non-dominated front 1 ($N_1$), \tiny Non-dominated front 0 ($N_0$)};
	\end{axis}
\end{tikzpicture}
 \caption{Non-dominated fronts}
    \label{fronts}
\end{figure}
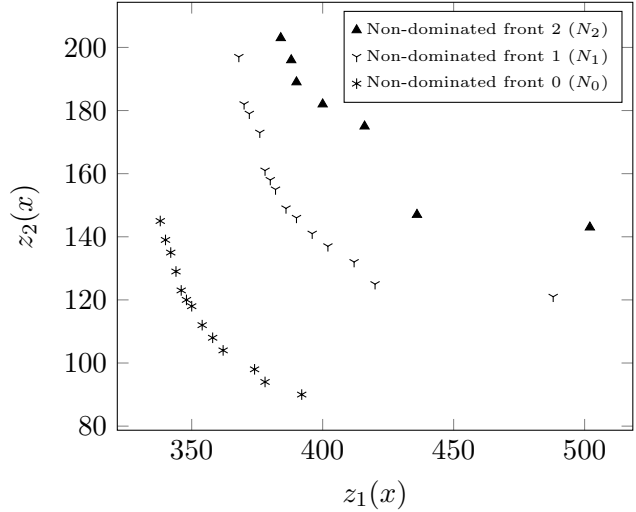

Before introducing the definition of crowding distance tournaments we still need to introduce the concept of successive non-dominated fronts. Let us denote here by $N_0$ the set of non-dominated solutions, or the non-dominated front. The remaining $n_f$ non-dominated fronts, $N_n,\; n=1,\ldots, n_n$, are the set of solutions that would be non-dominated if the solutions from the non-dominated fronts of smaller order are not considered. This concept can be visualised for a minimisation problem in Figure \ref{fronts}.

The crowding distance tournament is defined in Algorithm \ref{alg_cdt}. Its inputs are a set of solutions, $S$, and the number of rounds, $n_r$, and the output is the set of winners, $S^{win}$. In each round, a set of solutions, $S^{rand}$, is obtained by a random shuffle of the input set of solutions, $S$. The $shuffle()$ function is used to randomly modify the order of a set. It receives as input a set. Its output is a set with the same elements in a different order. Consecutive solutions, $S^{rand}_m,S^{rand}_{m+1}$, in the new list, $S^{rand}$, are paired for a match, $match(S^{rand}_m,S^{rand}_{m+1},S)$. The winner of this match is added to the set of winners, $S^{win}$. The match function is defined in Algorithm \ref{alg_match}.

\begin{algorithm}[h!]
\caption{Crowding distance tournament}\label{alg_cdt}
	\begin{algorithmic}[1]\normalsize
	\State{\textbf{input:} $S$ (list of solutions), $n_r$ (number of rounds)};
	\State{\textbf{output:} $S^{win}$ (winners)};
         \State{$S^{win}\leftarrow \{\}$;}
         \For{($r=1;r\leqslant n_r; +1$)}
         \State{$S^{rand}\leftarrow shuffle(S)$;}   
          \For{($m=1;m < length(S^{rand}); +2$)}
          \If{($length(S^{win})>0$)}
          \State{$S^{win}\leftarrow S^{win} \cup match(S^{rand}_m,S^{rand}_{m+1},S^{win})$;} 
          \Else
           \State{$S^{win}\leftarrow S^{win} \cup match(S^{rand}_m,S^{rand}_{m+1},S)$;} 
           \EndIf
        \EndFor
         \EndFor

        \State{{\textbf{return}($S^{win}$);}}
        \end{algorithmic} 
	   
\end{algorithm}

Let us denote by $N^{sol}$ a set where each element is the non-dominated front, $n$, where a solution, $s$, is located. There are two ways to win a match. The first is that one of the solutions is from a front with a lower rank than the other solution. If both solutions are from the same front, then the solution with the greater crowding distance to a certain set of solutions wins. If there are winners from previous matches, the considered set is the set of winners, $S^{win}$. Otherwise, the set of solutions, $S$, is considered. If there is a draw the first solution is chosen as the winner.

\begin{algorithm}[h!]
\caption{Crowding distance match}\label{alg_match}
	\begin{algorithmic}[1]\normalsize
	\State{\textbf{input:} $s^{1}$ (solution 1), $s^2$ (solution 2), $S$ (list of solutions)};
	\State{\textbf{output:} $s^{win}$ (winner)};
            \State{$s^{win}\leftarrow s^1$;} 
           \If{($N^{sol}_{s^2}<N^{sol}_{s^1}$)}
          \State{$s^{win}\leftarrow s^2$;} 
           \Else
           \If{($c(s^2,S)>c(s^1,S) \land N^{sol}_{s^2}=N^{sol}_{s^1}$ )}
          \State{$s^{win}\leftarrow s^2$;} 
           \EndIf
        \EndIf

        \State{{\textbf{return}($s^{win}$);}}
        \end{algorithmic} 
	   
\end{algorithm}
\subsection{NSGA-II}\label{nsga2}

\noindent The inputs of the NSGA-II algorithm are the problem being considered, $P$, the set of initial solutions, $S^{init}$, the population size, $n_{s}$, the number of generations, $n_g$, and the mutation probability, $m$. The output is a set of solutions to the problem, $S$. The pseudo-code for this algorithm is provided in Algorithm \ref{alg_2}. The first step of the algorithm is to initialise the population using Algorithm \ref{alg_i}. Then, for every generation, the offspring population, $S^{off}$, becomes the new parent population, $S^{par}$, from which a new offspring is generated. 
In each generation, the parent population is evaluated to determine the elite individuals $S^{elite}$. These individuals are automatically included in the new offspring. The procedure to find this set of individuals, $elite2(S^{par})$, is defined in Algorithm \ref{alg_e2}. The remaining individuals are generated through crossover and mutation. The parents are randomly chosen from a mating pool, $S^{pool}$. To compute which individuals from the elite population enter the mating pool a crowding distance tournament with 2 rounds takes place, $tournament(S^{elite},2)$. This procedure is defined in Algorithm \ref{alg_cdt}.

\begin{algorithm}[h!]
\caption{NSGA-II}\label{alg_2}
	\begin{algorithmic}[1]\normalsize
	\State{\textbf{input:} $P$ (problem), $S^{init}$ (initial solutions), $n_{s}$ (population size), $n_g$ (number of generations), $m$ (mutation probability)};
	\State{\textbf{output:} $S$ (set of solutions)};
        
        \State{$S^{off}\leftarrow initialise(n_{s},S^{init},A_{jt},B_i)$};
       \State{$S\leftarrow \{\}$;}

        \For{($g=1;g\leqslant n_g; +1$)}
         \State{$S^{par}\leftarrow S^{off}$};
        \State{$S^{elite}\leftarrow elite2(S^{par}, n_s)$};
        \State{$S^{pool}\leftarrow tournament(S^{elite},2)$};
        \State{$S^{off}\leftarrow S^{elite}$};
        \While{$length(S^{off})<n_{s}$}
        \State{$s^{p_1}\leftarrow rand(S^{pool}$)};
        \State{$s^{p_2}\leftarrow rand(S^{pool}$)};
        \State{$s^{off}\leftarrow crossover(s^{p_1}, s^{p_2})$};
        \State{$s^{off}\leftarrow mutation(s^{off}, m)$};
        \State{$S^{off}\leftarrow S^{off}\cup \{s^{off}$\}};
        \EndWhile
        \EndFor
         \State{$S\leftarrow S^{off}$};

        \State{{\textbf{return}($\mbox{$S$}$);}}
        \end{algorithmic} 
	   
\end{algorithm}

The inputs for finding the elite individuals of the population are the parent population, $S^{par}$, and the population size, $n_{s}$. Its output is the elite individuals, $S^{elite}$. The first step is to find the front where each individual belongs, as described in the previous subsection. Then, by ascending order of front, $n$, all of the individuals of that front, $N_n$, are added to the elite, until including an extra front would mean that the size of the elite would surpass half of the population size, $n_{s}/2$. Then, until the elite reaches this size, the individual from the first non-included front with the maximum crowding distance to the current elite is added.

\begin{algorithm}[h!]
\caption{NSGA-II elite}\label{alg_e2}
	\begin{algorithmic}[1]\normalsize
	\State{\textbf{input:} $S^{par}$ (parent population), $n_{s}$ (population size)};
	\State{\textbf{output:} $S^{elite}$ (elite individuals)};
       \State{$S^{elite}\leftarrow \{\}$};
       
       \State{$N\leftarrow evaluate\_fronts(S^{par})$};
         \State{$n' \leftarrow 1$};
       \For{($n = 0; n < length(N); +1$)}

       \If{($length(S^{elite} \cup N_n) \leqslant n_{s}/2$)}
       \State{$S^{elite} \leftarrow S^{elite} \cup N_n$};
       \Else
       \State{$n' \leftarrow n$};
       \State{\textbf{break}};
      
       \EndIf
    \EndFor

    \While{($length(S^{elite})<n_{s}/2$)}

     \State{$S^{elite} \leftarrow S^{elite} \cup \{\max crowding\_distance(N_{n'}, S^{elite})\}$};
    
    \EndWhile
    
    \State{{\textbf{return}($\mbox{$S^{elite}$}$);}}
    \end{algorithmic} 
	   
\end{algorithm}

After the elite is found, the mating pool, $S^{pool}$, is obtained through a 2-round crowding distance tournament (Algorithm \ref{alg_cdt}) between the elite individuals. Then, until the offspring has $n_{s}$ individuals, two parents, $s^{p_1}$ and $s^{p_2}$, are randomly selected from the mating pool. They undergo crossover to generate a new individual, $s^p$, through recombination of their genes (Algorithm \ref{alg_cross}). Then this individual is mutated with a probability, $m$, (Algorithm \ref{alg_mut}), and added to the offspring population, $S^{off}$. This whole process is repeated for each generation, until the last generation is found, and the algorithm returns the solutions of the last offspring population.

\subsection{NSGA-III}\label{nsga3}

\noindent Our implementation of the NSGA-III algorithm is similar to the NSGA-II implementation, with a key distinction in the definition of the elite of a population. This elite is defined based on the proximity of the individuals to a set of reference points, $Z^{ref}$. Algorithm \ref{alg_3} shows the pseudocode for the NSGA-III.

\begin{algorithm}[h!]
\caption{NSGA-III}\label{alg_3}
	\begin{algorithmic}[1]\normalsize
	\State{\textbf{input:} $P$ (problem), $S^{init}$ (initial solutions), $n_{s}$ (population size), $n_g$ (number of generations), $m$ (mutation probability), $Z^{ref}$ (reference points)};
	\State{\textbf{output:} $S$ (set of solutions)};
        
        \State{$S^{off}\leftarrow initialise(n_{s},S^{init},A_{jt},B_i)$};
       \State{$S\leftarrow \{\}$}

        \For{($g=1;g\leqslant n_g; +1$)}
         \State{$S^{par}\leftarrow S^{off}$};
        \State{$S^{elite}\leftarrow elite3(S^{par},Z^{ref})$};
        \State{$S^{pool}\leftarrow tournament(S^{elite},2)$};
        \State{$S^{off}\leftarrow S^{elite}$};
        \While{$length(S^{off})<n_{s}$}
        \State{$s^{p_1}\leftarrow rand(S^{pool}$)};
        \State{$s^{p_2}\leftarrow rand(S^{pool}$)};
        \State{$s^{off}\leftarrow crossover(s^{p_1}, s^{p_2})$};
        \State{$s^{off}\leftarrow mutation(s^{off}, m)$};
        \State{$S^{off}\leftarrow S^{off}\cup \{s^{off}$\}};
        \EndWhile
        \EndFor
         \State{$S\leftarrow S^{off}$};

        \State{{\textbf{return}($\mbox{$S$}$);}}
        \end{algorithmic} 
	   
\end{algorithm}

Algorithm \ref{alg_e3} provides the pseudocode for finding the elite in the NSGA-III algorithm. The initial step of adding non-dominated fronts to the elite until the first front, $n'$, that cannot be added without overfilling is the same as in the NSGA-II (Algorithm \ref{alg_e2}). However, in that algorithm the remaining spots would be filled based on the crowding distance of the solutions in the last front, $N_{n'}$, to the current elite. In this new algorithm, this is done based on the frequency of solutions assigned to a reference point. 

First, we find the set of solutions which should be assigned to each reference point, $Z^{assign}_{r'}$. These solutions come from the current elite individuals, $S^{elite}$, and the first non-included front, $N_{n'}$. To do this, we normalise the objective functions of each of these solutions based on the minimum and maximum values each objective takes. Then, each solution is assigned to the closest reference point. For example, if we have two reference points, (0.8, 0.2) and (0.2, 0.8), the solution (0.1, 0.6) would be assigned to the second reference point. After doing this assignment, we add individuals from the last front until the elite reaches the desired size. To select which individual to add, the first step is to evaluate the frequency of solutions already in the elite that are assigned to each of the points, $\omega_r$. Then, we find the point which has the least frequency of assignments, $r'$. Finally, we select a random point from the last front that has not yet been included and that is assigned to the reference point with the least frequency.

\begin{algorithm}[h!]
\caption{NSGA-III elite}\label{alg_e3}
	\begin{algorithmic}[1]\normalsize
	\State{\textbf{input:} $S^{par}$ (parent population), $n_{s}$ (population size), $Z^{ref}$ (reference points)};
	\State{\textbf{output:} $S^{elite}$ (elite individuals)};
       \State{$S^{elite}\leftarrow \{\}$};
   
       \State{$N\leftarrow evaluate\_fronts(S^{par})$};
         \State{$n' \leftarrow 1$};
       \For{($n = 0; n < length(N); +1$)}

       \If{($length(S^{elite} \cup N_n) \leqslant n_{s}/2$)}
       \State{$S^{elite} \leftarrow S^{elite} \cup N_n$};
       \Else
       \State{$n' \leftarrow n$};
       \State{\textbf{break}};
      
       \EndIf
    \EndFor
 \State{$Z^{assign}\leftarrow evaluate\_reference\_points(N_{n'} \cup S^{elite})$};
    \While{($length(S^{elite})<n_{s}/2$)}
     \State{$\omega_r\leftarrow evaluate\_frequency(S^{elite})$};
      \State{$r'\leftarrow \min \omega_r $};
     \State{$S^{elite} \leftarrow S^{elite} \cup\{ rand(Z^{assign}_{r'}\land N_{n'})\}$};
    
    \EndWhile
    
    \State{{\textbf{return}($\mbox{$S^{elite}$}$);}}
    \end{algorithmic} 
	   
\end{algorithm}
\section{Algorithmic framework}\label{s-AlgF}
\noindent This section introduces the decomposed multi-objective approach. It gives an overview of this method, compares it with other common methods, and explains its application to the single assignment thesis defence scheduling problem. Some required initiation protocols are addressed.

\subsection{Overview}\label{over}

\begin{figure}[h!]
  \centering

  \begin{subfigure}[b]{0.48\textwidth}
    \includegraphics[width=\textwidth]{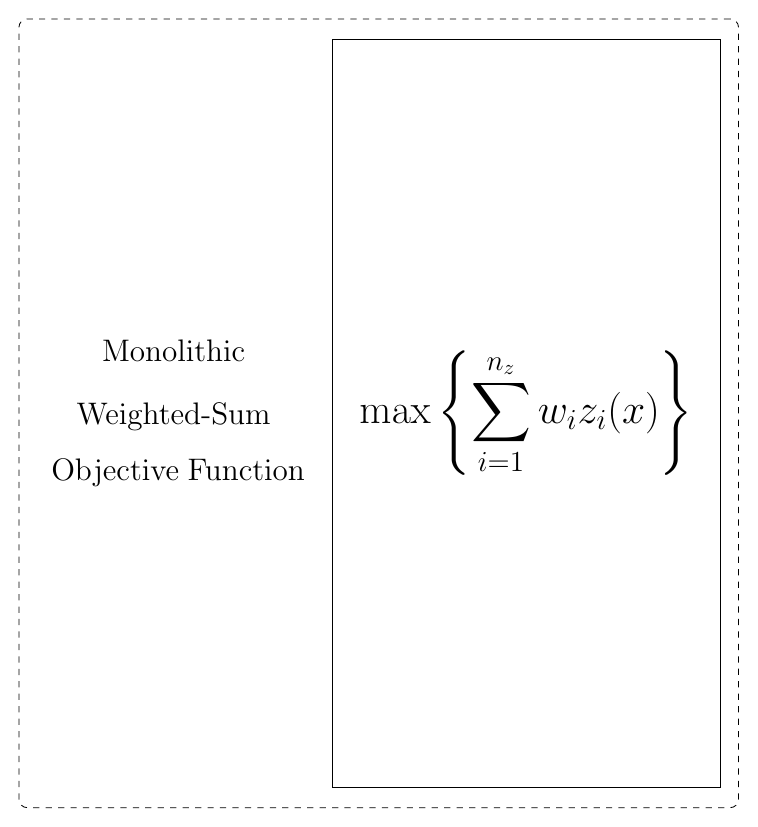}

  \end{subfigure}
  \hfill
  \begin{subfigure}[b]{0.48\textwidth}
    \includegraphics[width=\textwidth]{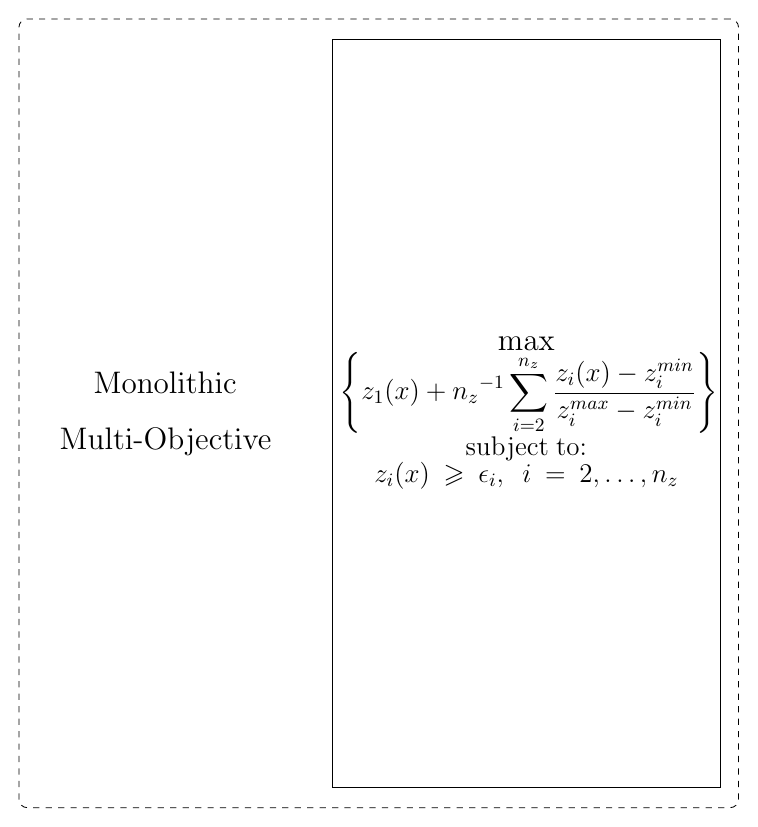}
 
  \end{subfigure}

  \vspace{1em}

  \begin{subfigure}[b]{0.48\textwidth}
    \includegraphics[width=\textwidth]{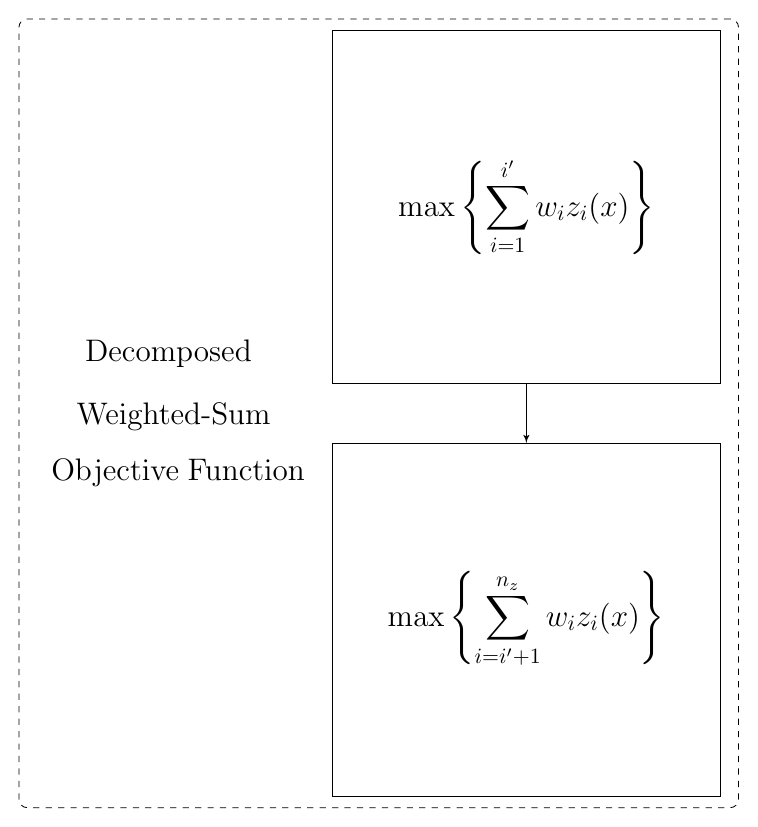}

  \end{subfigure}
  \hfill
  \begin{subfigure}[b]{0.48\textwidth}
    \includegraphics[width=\textwidth]{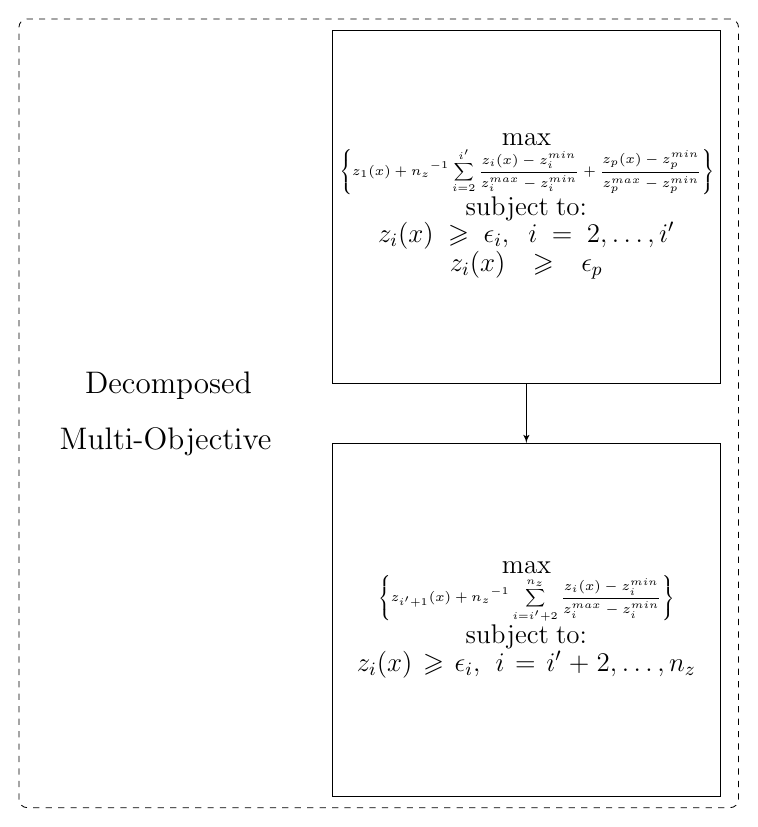}
 
  \end{subfigure}

  \caption{Method overview}
  \label{fig:main}
\end{figure}

\noindent Figure \ref{fig:main} displays four general methods to solve problems with multiple objectives. The monolithic weighted-sum objective function method considers a different weight, $w_i$, for each objective function, $z_i(x),\;i=1,\ldots,n_z$. Using this method only a single solution, which maximises this sum is provided to the decision-makers. The decomposed weighted-sum objective function works with similar objective functions. However, this type of method splits the problem into different stages to be solved sequentially. The problem is optimised for some objective functions, $z_i(x),\; i = 1,\ldots,i'$, and for some variables, $x \in X$, which define the performance of those objectives are fixed for the next optimisation stages, where the remaining objectives, $i = 1,\ldots,i'$, are optimised. This technique only guarantees the optimality of the objectives solved in the first stage. Nonetheless, for some problems, the solution time efficiency gained from the simplification of the models justifies the use of the method.

The multi-objective approaches have the advantage of providing decision-makers with a better knowledge of the non-dominated front but often come at the disadvantage of larger computational times. The monolithic multi-objective approaches, obtain non-dominated solutions for the entire problem considering all of the objectives.

In this work, we propose a decomposition approach for multi-objective optimisation. The aim is to improve its efficiency, while still finding good quality solutions. In this approach, the non-dominated solutions found in the first stage are used as partial solutions to be optimised in the second stage. This means that for every solution found in a previous stage, a new multi-objective algorithm is used to find the non-dominated solutions that can be derived from each partial solution. 

When decomposing problems with weighted-sum objective functions, the solutions are biased towards the objectives optimised in the first stage. In multi-objective decomposition, we propose the inclusion of proxy objectives, $z^{p}_j(x),\;j=1,\ldots,n_p$, that reflect the potential of a partial solution to have good performances in the objective functions of the next stage. This reduces the bias of the decomposition towards the initial stage objectives.

\subsection{Application}\label{app}

\noindent For a problem to be solved through this approach it should fulfil three assumptions. Like all multi-objective approaches, the problem must have several different objectives and conflicting points of view, making it advantageous for the decision-makers to be presented with multiple solutions with performance trade-offs. Similarly to the decomposed weighted-sum objective functions methods, the problem must be decomposable into sequential sub-problems where the previous stages provide partial solutions to the next. Finally, it must be possible to assess the partial solutions based on proxy objective functions that reflect the potential performance of the objectives that will be optimised in subsequent stages.

The single assignment thesis defence scheduling problem is a suitable candidate to showcase this method. It has multiple conflicting objectives, the problem can be split into sequential sub-problems, and the number of available slots for defences serves as a proxy to the objective functions of the second stage.

\begin{figure}[h!]
    \centering
    \includegraphics[scale=0.9]{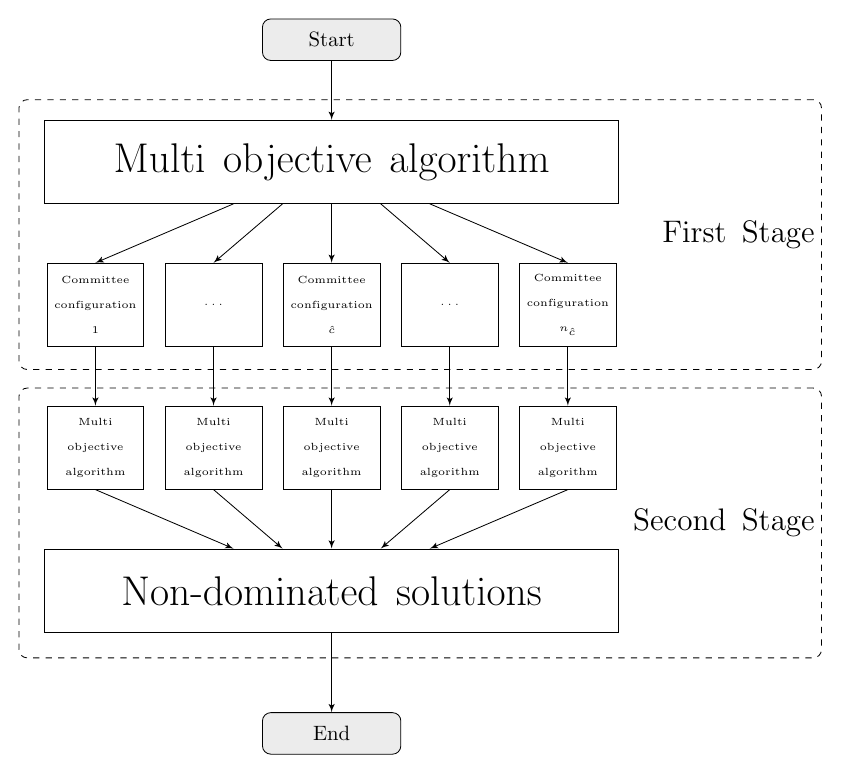}
    \caption{Decomposed multi-objective algorithms applied to thesis defence scheduling}
    \label{fig:enter-label}
\end{figure}

A diagram of this application is provided in Figure \ref{fig:enter-label}. In the first stage, a multi-objective algorithm (NSGA-II or NSGA-III), is used to find multiple committee configurations which are used as partial solutions in the next stage. In the second stage, for each configuration, a new multi-objective algorithm is used to find non-dominated time-slot assignments for each defence. Then all of the solutions found for each committee are assessed as a whole to find which are non-dominated considering the entire set and presented to the decision-makers. 

\subsection{Initialisation}\label{initi}

\noindent An initialisation step is required to find the minimum, $z^{min}$, and maximum, $z^{max}$, values vector for each of the considered objective functions before applying the $\epsilon$-constraint methods. For every objective, $z_i(x), \; i = 1,\ldots,n_z$, Objective Function \eqref{in_function} is used to obtain these values. The maximum values vector, $z^{max}$, takes the value of the respective objective function, $z^{max}_i=z_i(x), \;i'=1,\ldots,n_z$. The sum of the values of all other objectives is divided by a large enough value, $M$, such that this sum is smaller than 1 and does not interfere with the value of $z_{i}(x)$. For every objective, $z_i(x)$, the minimum values vector, $z^{min}_i$, takes the minimum value the objective, $z_i(x)$, took in this sum, $\min z_i'(x), \;i'=1,\ldots,n_z$. This ensures that the lower bounds are small enough such that the maximum values are attainable in the $\epsilon$-constraint method. In the genetic algorithms, the committee compositions that are obtained from each of these optimisations are included in the initial population.

\begin{equation}\label{in_function}
    \displaystyle \max z'_i(x) = z_i(x) + {M}^{-1}  \sum_{i'=1, i'\neq i}^{n_z} {z_i'(x)}.
\end{equation}
\section{Computational experiments}\label{s-exp}

\noindent This section presents the computational experiments. They aim to compare the performance of the decomposition strategy with a standard monolithic augmented $\epsilon$-constraint method. The experiments consider small and large size randomly generated instances, respectively. The small size instances have 25 committee members and 20 thesis defences. For each defence, 2 committee members are already assigned and the model selects 1 additional member. The large size instances have 50 committee members and 40 thesis defences. For each defence, one committee member is already assigned, and the model selects two additional members. Additional information regarding the generation of these instances is provided in \cite{Almeida2024}. 

Each augmented $\epsilon$-constraint iteration is conducted with a time limit of 120 seconds. The CPU is an Intel(R) Core(TM) i7-8565U CPU @ 1.80GHz   1.99 GHz, and the installed RAM is 8GB. The genetic algorithms are implemented in C++ and the augmented $\epsilon$-constraints are implemented in Python and solved using Gurobi.

\subsection{Small size randomly generated instances}\label{srng}

\noindent The experiments presented in this subsection concern Instances 5 and 6 from \cite{Almeida2024}. We test parameterisations for our method and compare their performance against a monolithic $\epsilon$-constraint method. Let us note that the initialisation phase is the same for each experiment. In this phase, four solutions are found by maximising the monolithic problem considering Objective Function \eqref{in_function}. 

\begin{figure}[H]
    \centering
    \includegraphics[scale=0.3]{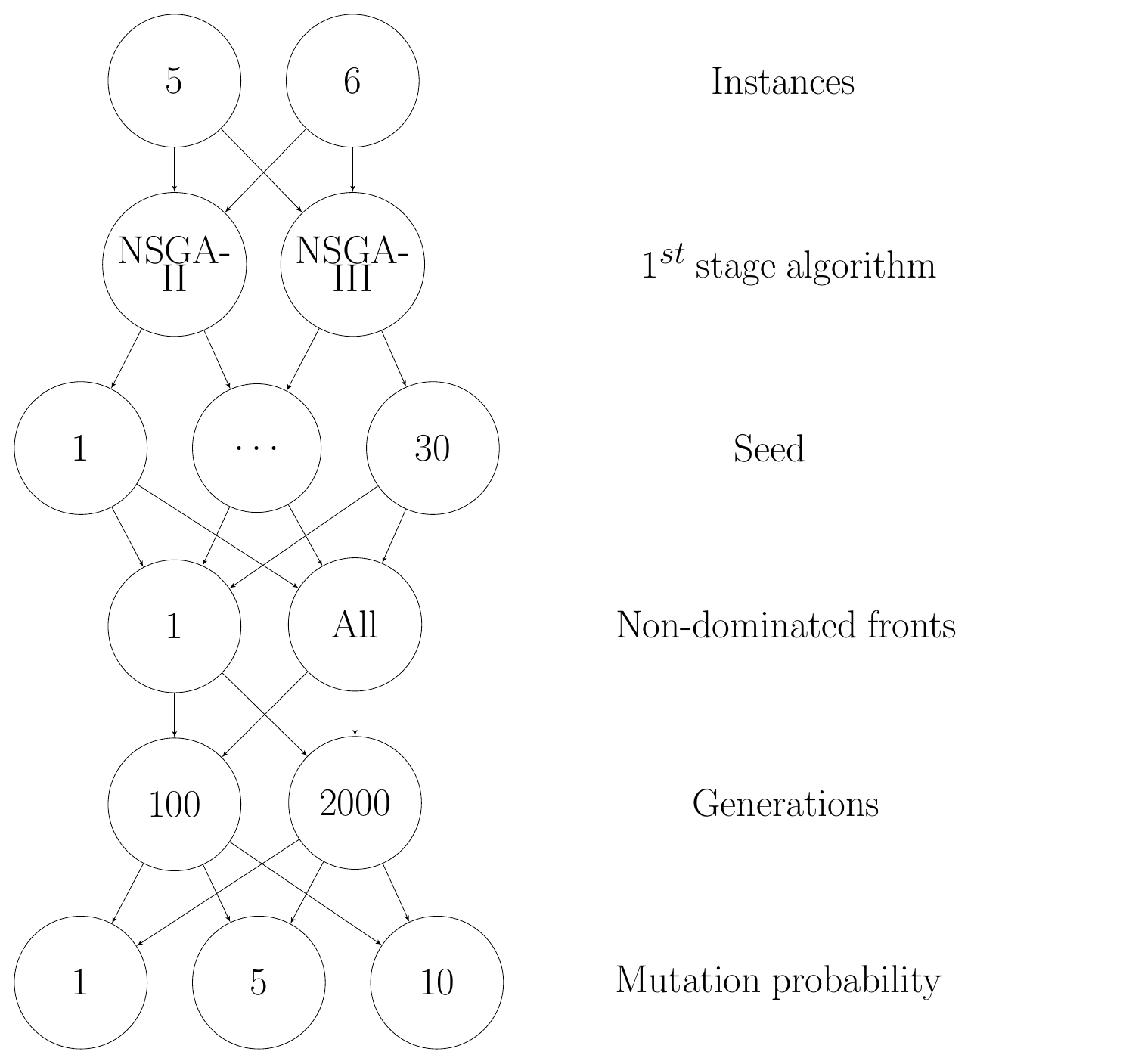}
    \caption{Computational experiments with small size randomly generated instances}
    \label{small}
\end{figure}

Instead of solving a monolithic problem, we decompose it and find partial solutions in the first stage using genetic algorithms. In the second stage, for each partial solution, non-dominated solutions are found using an augmented $\epsilon$-constraint method. Partial solutions are generated using different sets of generations, and mutation probabilities. The population size is always set to 200. We also consider using only the non-dominated front of partial solutions and using all the partial solutions that are found. For each parameterisation, 30 sets of partial solutions are generated with different seeds to ensure the statistical relevance of the results. A summary of the experiment design process for the generation of the sets of partial solutions is provided in Figure \ref{small}.

The results for Instances 5 and 6 are presented in Tables \ref{tab:performance_comparison} and \ref{tab:performance_comparison6}. For the monolithic resolutions, the tables include information regarding the computational time, hyper-volume, and number of non-dominated solutions found. For the decomposed resolutions, the tables include information regarding the average times of each stage, average hyper-volume, average number of solutions found that are not dominated by any solution found in the monolithic resolutions ($N^{mono}_0$), and the average number of non-dominated solutions found ($N_0$). Let us note that each of these rows pertains to 30 runs with different seeds and their average results. The hyper-volume for the point considering the minima of each objective and the initialisation time are also provided in the table notes.

\begin{table}[H]
  \centering
  \caption{Performance comparison of methods for Instance 5 with $n_{s}=200$ and $n_g=2000$}
  \label{tab:performance_comparison}
  \scalebox{0.95}{
   \begin{tabular}{|c|c|c|c|c|c|c|c|c|c|c|c|c|}
    \hline
    \multirow{2}{*}{Algorithm}& \multirow{2}{*}{$m$}& \multirow{2}{*}{$n_f$}&\multicolumn{2}{c|}{$1^{st}$ stage(sec)} & \multicolumn{2}{c|}{$2^{nd}$ stage(sec)} &\multicolumn{2}{c|}{Hyper-volume} & \multicolumn{2}{c|}{$N^{mono}_0$}&\multicolumn{2}{c|}{$N_0$} \\
 \cline{4-13}
   &&&100&2000&100&2000&100&2000&100&2000&100&2000\\ 
    \hline
    & 1     & 1     &   \multirow{6}{*}{7}               & \multirow{6}{*}{134}       &    209    &  212       &0.1155  &0.1164   &33& 53& 210 & 200 \\
      & 5     & 1     &                    &       &    220&       232            & 0.1157 &\textbf{0.1167} &41& 55  &213&164  \\
    NSGA-II+& 10    & 1     &                    &      & 213&   234              & \textbf{0.1159}   & \textbf{0.1167}  &39 &52&201&  162  \\
   \cdashline{2-3}\cdashline{6-13}
      $\epsilon$-constraint($\epsilon'=1)$ & 1     & All   &                    &       &291&310               &0.1156    & 0.1165 &34   &55 &215&207 \\
     & 5     & All   &                    &        &    294&341              &0.1158&\textbf{0.1168} &42&56&221&173  \\
      & 10    & All   &                    &        &302& 345             & \textbf{0.1161}    &\textbf{0.1168}  &40  &53&212&173  \\
     
     \hdashline
     {}   & 1     & 1     &      \multirow{6}{*}{7}              &  \multirow{6}{*}{135}   &   205&     217             & 0.1153 & 0.1161&38  & 40&297&  224  \\
     & 5     & 1     &                    &    &    202&       221            & 0.1159&\textbf{0.1167} &47&42 &242& 210  \\
 NSGA-III+  & 10    & 1     &                    &        &   200&  222          & \textbf{0.1160} & \textbf{0.1167}   &56  & 40&264&  208 \\
      \cdashline{2-3}\cdashline{6-13}
     $\epsilon$-constraint($\epsilon'=1)$ & 1     & All   &                    &        &   291&   300           &0.1154 & 0.1162 &39&41&  299 &229  \\
    & 5     & All   &                    &       &    301&315              &0.1160&\textbf{0.1168} &47&42& 245 &219  \\
      & 10    & All   &                    &        &    295&  329            & \textbf{0.1162}&\textbf{0.1168} &57&41& 273& 216  \\
    
    \hdashline
    Monolithic &
   \multicolumn{2}{c|}{\multirow{2}{*}{-}}     &     \multicolumn{4}{c|}{\multirow{2}{*}{1394}}&\multicolumn{2}{c|}{\multirow{2}{*}{0.1206}}            & \multicolumn{4}{c|}{\multirow{2}{*}{89}} \\
   ($\epsilon'=\frac{z^{max}-z^{min}}{10})$&\multicolumn{2}{c|}{}&\multicolumn{4}{c|}{}&\multicolumn{2}{c|}{}&\multicolumn{4}{c|}{}\\
    \hdashline
    Monolithic &
    \multicolumn{2}{c|}{\multirow{2}{*}{-}}     &     \multicolumn{4}{c|}{\multirow{2}{*}{72356 }}&\multicolumn{2}{c|}{\multirow{2}{*}{0.1214}}            & \multicolumn{4}{c|}{\multirow{2}{*}{1301}} \\
    ($\epsilon'=1)$&\multicolumn{2}{c|}{}&\multicolumn{4}{c|}{}&\multicolumn{2}{c|}{}&\multicolumn{4}{c|}{}\\
    \hline

  \end{tabular}}
  {\raggedright \footnotesize $^1$Results for the genetic algorithms are averages for 30 runs with different seeds\\
  $^2$Minimum hyper-volume=0.0420\\
  $^3$Initialisation time = 16 seconds\\}
\end{table}

\begin{table}[H]
  \centering
  \caption{Performance comparison of methods for Instance 6 with $n_{s}=200$ and $n_g=2000$}
  \label{tab:performance_comparison6}
  \scalebox{0.95}{
 \begin{tabular}{|c|c|c|c|c|c|c|c|c|c|c|c|c|}
    \hline
    \multirow{2}{*}{Algorithm} & \multirow{2}{*}{$m$} & \multirow{2}{*}{$n_f$} & \multicolumn{2}{c|}{{$1^{st}$ stage(sec)}} & \multicolumn{2}{c|}{{$2^{nd}$ stage(sec)}} & \multicolumn{2}{c|}{{Hyper-volume}} & \multicolumn{2}{c|}{{$N^{mono}_0$}} & \multicolumn{2}{c|}{{$N_0$}} \\
    \cline{4-13}
    &&& 100 & 2000 & 100 & 2000 & 100 & 2000 & 100 & 2000 & 100 & 2000 \\
    
    \hline
     {}   &1 &1&\multirow{6}{*}{7}&\multirow{6}{*}{134}&244&  239& 0.1247   &0.1254&44&47& 268&  260\\ 
    &5&1&& &   239&    263& 0.1252 & 0.1263 &46  &41&275&  228\\
  NSGA-II+&10&1&& &    236&  270 & \textbf{0.1253}&  \textbf{0.1269} &45&32&261& 213 \\ 
   
\cdashline{2-3}\cdashline{6-13}
  $\epsilon$-constraint($\epsilon'=1)$&1&All&& & 322& 348&0.1249&0.1255&46&49&288&  279   \\
   &5&All&&  &    351&    382&0.1253&0.1263&47&42&290&  249 \\
   &10&All&& & 351&   381&\textbf{0.1254}& \textbf{0.1269} &45& 33&281&  235\\
   
\hdashline
&1&1&\multirow{6}{*}{7}&\multirow{6}{*}{137}&    239&    228 &0.1249 & 0.1253&51&54 &311& 341   \\
   &5&1&& &   240&  235&\textbf{0.1250}&0.1259&52&47 &310& 320  \\
  NSGA-III+&10&1&&  &    218&    242&\textbf{0.1250}&\textbf{0.1263}&46&45&294& 285 \\
 \cdashline{2-3}\cdashline{6-13}
   $\epsilon$-constraint($\epsilon'=1)$&1&All&& &  344
& 334&0.1249& 0.1254&52&55&329& 358   \\
 &5&All&& &343&   358&\textbf{0.1251}& 0.1260&53&47 &332& 347 \\
  &10 &All&& & 313&  361&\textbf{0.1251}&\textbf{0.1264}&46&45&314&  309  \\

    \hdashline
    
   Monolithic & \multicolumn{2}{c|}{\multirow{2}{*}{-}}     &     \multicolumn{4}{c|}{\multirow{2}{*}{3162}}&\multicolumn{2}{c|}{\multirow{2}{*}{0.1301  }}            & \multicolumn{4}{c|}{\multirow{2}{*}{186}} \\
   ($\epsilon'=\frac{z^{max}-z^{min}}{10})$&\multicolumn{2}{c|}{}&\multicolumn{4}{c|}{}&\multicolumn{2}{c|}{}&\multicolumn{4}{c|}{}\\
    \hline

  \end{tabular}}
  {\raggedright \footnotesize $^1$Results for the genetic algorithms are averages for 30 runs with different seeds\\
  $^2$Minimum hyper-volume=0.0103\\
  $^3$Initialisation time = 39\\}
\end{table}

The decomposition method is considerably faster than the monolithic resolutions. The experiments with decomposition took between 8\%  and 32\% of the time of monolithic resolutions considering 10 increments between each $z^{min}$ and $z^{max}$. These time differences are due to differences between instances, the number of generations in the first stage, which naturally affects its running time, and the number of non-dominated fronts of partial solutions being considered, which impact the number of augmented $\epsilon$-constraint iterations for the second stage. There is a trade-off between this time and the hyper-volume of the objective space. Nonetheless, the solutions found by using decomposition are competitive with the ones found using the monolithic approach. 

Allowing the genetic algorithms to run for more generations led to higher hyper-volume values in both instances. This is an expected outcome, but it also indicates that the improvement of the proxy objective over the generations leads to better solutions being found. Higher mutation probabilities also promoted higher hyper-volume values. This is caused by the increased diversification of solutions. For Instance 5, NSGA-III provided better hypervolumes and for Instance 6, NSGA-II had the better performance.

Decomposition can be useful for small size multi-objective problems if finding solutions faster is important for the decision-makers. Nonetheless, this comes at the cost of having slightly worse solutions to choose from than when using a more time-consuming method such as a monolithic augmented $\epsilon$-constraint.

\subsection{Large size randomly generated instances}\label{lrng}

\begin{figure}[h]
    \centering
    \includegraphics[scale=0.3]{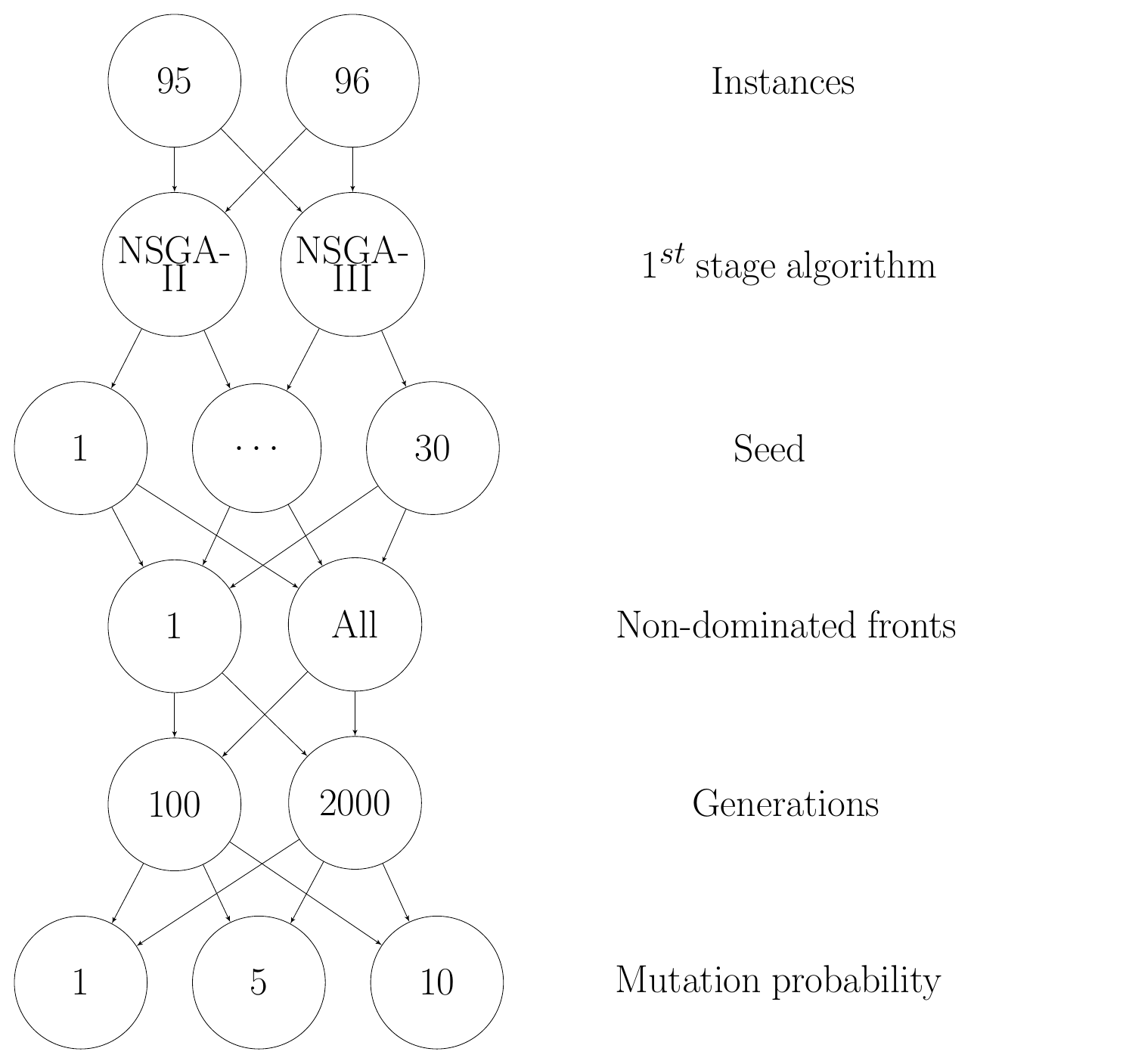}
    \caption{Computational experiments with large size randomly generated instances}
    \label{large}
\end{figure}

\noindent The experiments presented in this subsection consider Instances 95 and 96 from \cite{Almeida2024}. These instances have more committee members and thesis defences than the ones from the previous subsection. There are also more roles to be assigned to complete each committee. The tested parameterisations (Figure \ref{large}) and result presentation (Tables \ref{tab:performance_comparison95} and \ref{tab:performance_comparison96}) follow the same structure as in the previous subsection.

\begin{table}[H]
  \centering
  \caption{Performance comparison of methods for Instance 95 with population size $n_{s}=200$, and number of generations $n_g=2000$}
  \label{tab:performance_comparison95}
   \scalebox{0.95}{
  \begin{tabular}{|c|c|c|c|c|c|c|c|c|c|c|c|c|}
    \hline
    \multirow{2}{*}{Algorithm} & \multirow{2}{*}{$m$} & \multirow{2}{*}{$n_f$} & \multicolumn{2}{c|}{{$1^{st}$ stage(sec)}} & \multicolumn{2}{c|}{{$2^{nd}$ stage(sec)}} & \multicolumn{2}{c|}{{Hyper-volume}} & \multicolumn{2}{c|}{{$N^{mono}_0$}} & \multicolumn{2}{c|}{{$N_0$}} \\
    \cline{4-13}
    &&& 100 & 2000 & 100 & 2000 & 100 & 2000 & 100 & 2000 & 100 & 2000 \\
    \hline
   &  
    1 & 1 &\multirow{6}{*}{16}  & \multirow{6}{*}{313} &774& 758 &0.1501& \textbf{0.1507} &159& {113} &491& {451} \\
    & 5 & 1 &  &  &769& 810 &\textbf{0.1505}& 0.1500 &151& 93 &507& 442 \\
    NSGA-II+& 10 & 1 &  &  &756& 754 &0.1503& 0.1500 &139& 86 &484& 380 \\
  \cdashline{2-3}\cdashline{6-13}
     $\epsilon$-constraint($\epsilon'=1)$ & 1 & All &  &  &1078& 1014 &0.1504& \textbf{0.1509} &164& {117} &557& {513} \\
    & 5 & All &  &  &   1024& 1085 & \textbf{0.1507} & 0.1503 &155 & 96 &568& 498 \\
    & 10 & All &  &  &  1055& 1061 & 0.1506  & 0.1502 &146  & 87 &553& 428 \\
    
    \hdashline
     & 
    1 & 1 &\multirow{6}{*}{17}  & \multirow{6}{*}{319} &605& 578 &0.1492& 0.1494 &191& 123 &448& 365 \\
    & 5 & 1 &  &  &598& 617 &0.1493& \textbf{0.1502} &190& {147} &439& 414 \\
   NSGA-III+ & 10 & 1 &  &  &607& 689 &\textbf{0.1497}& 0.1500 &193& 131 &450& {456} \\
   \cdashline{2-3}\cdashline{6-13}
   $\epsilon$-constraint($\epsilon'=1)$ & 1 & All &  &  &   885
& 799 & 0.1496 & 0.1499 &199& 128 &510& 425 \\
    & 5 & All &  &  & 856& 857 &0.1497 & \textbf{0.1507} &194 & {154} &499& 474 \\
    & 10 & All &  &  &874& 974 &\textbf{0.1500}    & 0.1503 &202  & 134 &517 & {503} \\
   
    \hdashline
    
     Monolithic & \multicolumn{2}{c|}{\multirow{2}{*}{-}}     &     \multicolumn{4}{c|}{\multirow{2}{*}{9200}}&\multicolumn{2}{c|}{\multirow{2}{*}{0.1475}}            & \multicolumn{4}{c|}{\multirow{2}{*}{55}} \\
   ($\epsilon'=\frac{z^{max}-z^{min}}{10})$&\multicolumn{2}{c|}{}&\multicolumn{4}{c|}{}&\multicolumn{2}{c|}{}&\multicolumn{4}{c|}{}\\
    \hline
  \end{tabular}}
  {\raggedright \footnotesize $^1$Results for the genetic algorithms are averages for 30 runs with different seeds\\
  $^2$Minimum hyper-volume=0.0153\\
  $^3$Initialisation time = 293\\}
\end{table}

\begin{table}[H]
  \centering
  \caption{Performance comparison of methods for Instance 96 with population size $n_{s}=200$, and number of generations $n_g=100$ and $n_g=2000$}
  \label{tab:performance_comparison96}
   \scalebox{0.95}{
  \begin{tabular}{|c|c|c|c|c|c|c|c|c|c|c|c|c|}
    \hline
    \multirow{2}{*}{Algorithm}& \multirow{2}{*}{$m$}& \multirow{2}{*}{$n_f$}&\multicolumn{2}{c|}{$1^{st}$ stage(sec)} & \multicolumn{2}{c|}{$2^{nd}$ stage(sec)} &\multicolumn{2}{c|}{Hyper-volume} & \multicolumn{2}{c|}{$N^{mono}_0$}&\multicolumn{2}{c|}{$N_0$} \\
 \cline{4-13}
   &&&100&2000&100&2000&100&2000&100&2000&100&2000\\ 
    \hline
    &  
    1&	1&	\multirow{6}{*}{17}&				\multirow{6}{*}{314}	&568	&611	&\textbf{0.1963}	&\textbf{0.1961}	&{143}	&{147}	&{314}	&{324}	\\
&	5	&	1	&	&	&536&657	&0.1958	&0.1950	&137	&116	& 291	&315	\\
 {NSGA-II+}&	10	&	1	&	&	&543&698	&0.1959	&0.1945	&	125&110	&286	&321	\\
\cdashline{2-3}\cdashline{6-13}
{$\epsilon$-constraint($\epsilon'=1)$}&	1	&	All	&	&	&868&862	&\textbf{0.1967}	&\textbf{0.1965}	&{168}	&{167}	&{360}	&{365}	\\

&	5	&	All	&	&	&810&912	&0.1963	&0.1955	& 159	&133	&340	&352	\\

&	10	&	All	&	&	&849&985	&0.1964    	&0.1950	&149  	&122	&331 	&357	\\

   \hdashline
      &1&	1&\multirow{6}{*}{17}	&\multirow{6}{*}{318}&   543&565	&\textbf{0.1952}	&\textbf{0.1956}	&172   	&{197}	&{315}	&349	\\
&	5	&	1	&	&	&   502&579	& \textbf{0.1952} 	&0.1942	&{177} &182	&312		&{363}	\\
NSGA-III+&10&1&&&	 511&596	&0.1949	&0.1933	&160 &150	& 301 		&340	\\
\cdashline{2-3}\cdashline{6-13}
$\epsilon$-constraint($\epsilon'=1)$&	1	&	All	&	&	&   777&793	&\textbf{0.1958}	&\textbf{0.1962}	&	205   &{224}	&{365}	&392	\\
&	5	&	All	&	&	&    765&838	& 0.1957	&0.1948	&{209} 	&206	&364	&{404}	\\
&	10	&	All	&	&	& 786&886	& 0.1955   	&0.1938	&	189  &170	&350	&379	\\

     \hdashline
    Monolithic & \multicolumn{2}{c|}{\multirow{2}{*}{-}}     &     \multicolumn{4}{c|}{\multirow{2}{*}{13380}}&\multicolumn{2}{c|}{\multirow{2}{*}{0.1926}}            & \multicolumn{4}{c|}{\multirow{2}{*}{82}} \\
   ($\epsilon'=\frac{z^{max}-z^{min}}{10})$&\multicolumn{2}{c|}{}&\multicolumn{4}{c|}{}&\multicolumn{2}{c|}{}&\multicolumn{4}{c|}{}\\
    \hline

  \end{tabular}}
  {\raggedright \footnotesize $^1$Results for the genetic algorithms are averages for 30 runs with different seeds\\
  $^2$Minimum hyper-volume=0.0192\\
  $^3$Initialisation time = 280 seconds\\}
\end{table}

The running time for the experiments with decomposition was between 6\% and 18\% of the running time for the monolithic experiments. In contrast with the less complex instances from the previous subsection, the hyper-volume values found while decomposing the problem were better than when solving the monolithic version. This happens because in these instances the solver is not always able to reach a 0\% gap in each monolithic iteration due to the 120 seconds time limit. If this time limit was not imposed, the monolithic method would find better results, but this would be too impractical due to the large computational time it would require.

Similarly to the smaller size instances, Instance 95 benefited from longer first stages, even if the difference was not as large. However, the same was not true for Instance 96. In this case, some sets of partial solutions found after 100 generations yielded better final solutions than those found after 2000 generations. Due to the increased instance complexity, some relations and interactions between different assignments might not be as well represented by the proxy objective as for the simpler instances. This means that parameterisations leading to partial solutions which are more distinct from the ones found using an exact approach in the initialisation phase end up losing some of these interactions and producing worse solutions. This is supported by the fact that in both of these instances, but especially for Instance 96, smaller mutation probabilities lead to higher hyper-volume values.

Decomposition seems especially well suited for larger size instances where monolithic resolutions take too much time and can end up with worse solutions. However, finding the right parameterisations is a complex process, as there is a trade-off between allowing the first stages to run for longer periods and improving the objective functions that are assessed in this stage, and losing some of the interactions that allow for better performance of the objectives that are only set in the second stage.

To sum up, we can derive the following conclusions from these computational experiments:
\begin{itemize}[label={--}]
\item Taking a decomposition approach can considerably improve the efficiency of multi-objective methods.
\item There is a trade-off between this increased efficiency and the quality of solutions. This approach is better suited for more complex instances where monolithic approaches fail to reach optimality.
\item In most cases allowing a longer first stage improves the quality of solutions.
\item Certain interdependencies and interactions are not well represented by the proxy objective. This means that shorter first stages can produce better results for some instances.
\end{itemize}
\section{Real-world case study}\label{s-cases}

\noindent This section analyses a real-world case study. It is based on the thesis scheduling process of the Department of Engineering and Management at the Instituto Superior Técnico of the Universidade de Lisboa. The instance, data gathering, and some algorithmic adaptations are presented. Before analysing the results, an explanation of how to interpret the objective function values is provided.

\subsection{Instance and data gathering}\label{idg}

\noindent The case study has 47 committee members, 36 thesis defences, 3 roles, the supervisor, whom is already assigned, the president and a third member. There are 16 different days, with 31 time slots (15 minutes \textit{per} slot), and 2 rooms. The duration of each defence is one hour, or 4 time slots. Three objectives are considered, balanced workloads, $z_1(x)$, time slot preferences, $z_3(x)$, and scheduled days, $z_4(x)$.

When the committee member availability data was gathered by the person responsible for scheduling the defences the committees for each defence were already assigned. For each defence, doodles were sent to the supervisor and the third member with slots where the president was available and each member selected those they were also available for. For some defences, this method did not yield any agreeable slot, and the person responsible had to try to mediate a solution that was acceptable for all three members. This is a considerably taxing procedure for all the people involved, and the results are not always great.

\subsection{Algorithmic adaptations}\label{aa}

\noindent When doing the initial tests for this instance we realised that some particular interactions and relations between committee assignments (partial solutions) had a more significant impact on the quality of the final schedules when compared to the randomly generated instances analysed in the previous section. This is perhaps due to the data-gathering process, where the availability of committee members is assessed based on pre-established committee assignments. This means that many assignments that would be possible are not feasible considering the availability data that we have access to. This limitation means that the more random initial population generation method did not yield very good results for this particular instance.

To ensure that these interdependencies are represented in the initial population, the generation of the initial population is based on crossovers between the solutions found in the initialisation phase by optimising the monolithic problem for each of the objective functions. This crossover operator is presented in Algorithm \ref{alg_nucross}. It is similar to the one previously presented in Algorithm \ref{alg_cross}. However, the probability of a parent chromosome being chosen for each committee is an input, instead of being set at 50\%. 

\begin{algorithm}[h!]
\caption{Non-uniform crossover}\label{alg_nucross}
	\begin{algorithmic}[1]\normalsize
	\State{\textbf{input:}  $s^{p_1}$ (parent 1), $s^{p_2}$ (parent 2), $v$ (probability)};
	\State{\textbf{output:} $s^{off}$ (offspring)};

        \State{$s^{off}\leftarrow\{\}$;}
        \For{($j=1;j\leqslant n_j; +1$)}
        \State{$s^p\leftarrow rand(s^{p_1},s^{p_2},v)$;}
         \State{$s^{off}\leftarrow s^{off} \frown s^p_j$;}
        \EndFor

        \State{{\textbf{return}($\mbox{$s^{off}$}$);}}
        \end{algorithmic} 
	   
\end{algorithm}

New individuals are generated based on each combination of initial solutions ($s^i$, $s^j$). A set of probabilities, $V$, is defined and $n_h$ new individuals are generated based on non-uniform crossovers for each probability, $v$, within that set. This process can be compared to a path-relinking procedure between the partial solutions, ($s^i$, $s^j$).

\begin{algorithm}[h!]
\caption{Adapted initialisation}\label{alg_ai}
	\begin{algorithmic}[1]\normalsize
	\State{\textbf{input:}  $S^{init'}$ (initial solutions), $n_r$ (repetitions), $V$ (probability vector)};
	\State{\textbf{output:} $S^{init}$ (initial solutions/population)};

        \State{$S^{init}\leftarrow S^{init'}$;}
        \For{($s^i\in S^{init'}$)}
        \For{($j=i+1;j\leqslant n_i; +1$)}
        \For{($v\in V$)}
        \For{($h=1;r\leqslant n_h;+1$)}
         \State{$S^{init}\leftarrow S^{init} \cup \{non\_uniform\_crossover(s^i,s^j,v)\}$;}
        \EndFor
         \EndFor
          \EndFor
         \EndFor      
        \State{{\textbf{return}($\mbox{$S^{init}$}$);}}
        \end{algorithmic} 
	   
\end{algorithm}
\subsection{Objective function interpretation}\label{ofi}

\noindent This subsection aims to help the reader in interpreting the meaning of the different objective values to better understand the implications of the results presented in the following subsection. The quality of the solutions in this case study is assessed based on three objectives, the workload balance objective, $z_1(x)$, the time slot preference objective, $z_3(x)$, and the committee days objective, $z_4(x)$. The preference slots objective, $z_3(x)$, is straightforward to interpret. The value of the objective represents the number of undesirable time slots where a member has a defence scheduled. For example, a value of $z_3(x)=-8$, might mean that 2 members have a defence scheduled in undesirable time slots. Let us note that each defence lasts for 4 time slots.

\begin{figure}[H]
        \centering
    \scalebox{0.8}{
   \hspace{-0cm} \begin{tikzpicture}

\pie[pos={10,0}, radius=2.5, color={red!10, red!20, red!30, red!40, red!50,red!60,red!70,red!80}]
{75/, 13/, 13/}

\pie[pos={16,0}, radius=2.5, color={red!5,red!10, red!20, red!30, red!40, red!50,red!60,red!70,red!80}]
{28/, 34/, 17/, 17/, 4/}

\node[block3, yshift=3cm, xshift= 10cm]{$z_1(x)=-338$};

\node[block3, yshift=3cm, xshift= 16cm]{$z_1(x)=-488$};

\node[block3, yshift=2cm, xshift= 20.5cm]{Assignments};
\node[legend0, xshift=20.5cm,yshift=1.2cm](0){0};
\node[legend1, xshift=20.5cm,yshift=0.6cm](1){1-2};
\node[legend2, xshift=20.5cm,yshift=0cm](2){3-4};
\node[legend3, xshift=20.5cm,yshift=-0.6cm](1){5-6};
\node[legend4, xshift=20.5cm,yshift=-1.2cm](1){7-8};

\end{tikzpicture}}
    \caption{Percentage of committee members with a number of committees in two solutions with different workload balance objective values}
    \label{Z1 deg}
\end{figure}
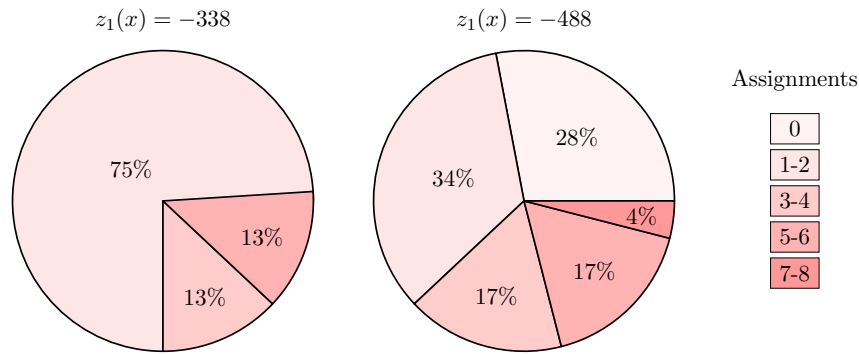

The balanced workloads objective, $z_1(x)$, is a measure of fairness regarding the number of assignments for committee members. Higher absolute values imply more uneven numbers of assignments. Figure \ref{Z1 deg} compares the number of assignments \textit{per} member in two solutions with different performances in this objective. In the solution with better performance, $z_1(x)=-338$, most members have 1 or 2 assignments. In the other solution, $z_1(x)=-488$, a higher percentage of members have more assignments, 38\% of members have 3 or more assignments, whereas in the other solution, this percentage was 25\%. This solution also has some members without assignments and some members with 7 or 8 assignments.

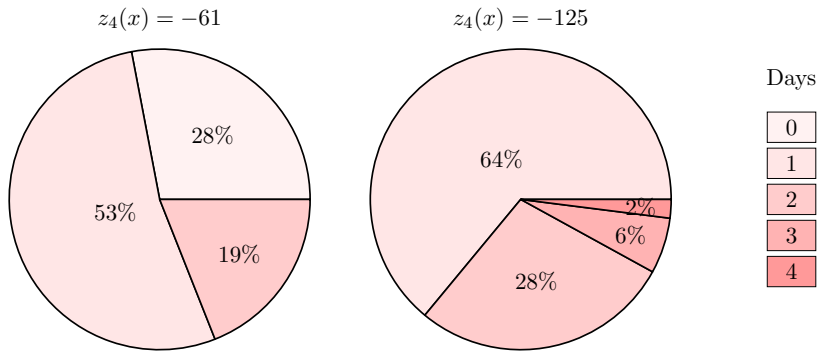
\begin{figure}[H]
        \centering
   \scalebox{0.8}{ \begin{tikzpicture}

\pie[pos={10,0}, radius=2.5, color={red!5, red!10, red!20, red!30, red!40, red!50,red!60,red!70,red!80}]
{28/, 53/, 19/}

\pie[pos={16,0}, radius=2.5, color={red!10, red!20, red!30, red!40, red!50,red!60,red!70,red!80}]
{64/, 28/, 6/, 2/}

\node[block3, yshift=3cm, xshift= 10cm]{$z_4(x)=-61$};

\node[block3, yshift=3cm, xshift= 16cm]{$z_4(x)=-125$};
\node[block3, yshift=2cm, xshift= 20.5cm]{Days};
\node[legend0, xshift=20.5cm,yshift=1.2cm](0){0};
\node[legend1, xshift=20.5cm,yshift=0.6cm](1){1};
\node[legend2, xshift=20.5cm,yshift=0cm](2){2};
\node[legend3, xshift=20.5cm,yshift=-0.6cm](1){3};
\node[legend4, xshift=20.5cm,yshift=-1.2cm](1){4};

\end{tikzpicture}}
    \caption{Percentage of committee members with a number of scheduled days in two solutions with different committee days objective values}
    \label{Z4 deg}
\end{figure}

The committee days objective, $z_4(x)$, aims to fairly minimise the number of days committee members have defences scheduled on. Higher absolute values imply more assigned days. Figure \ref{Z4 deg} compares the number of assignments \textit{per} member in two solutions with different performances in this objective. In the solution with the better performance, $z_4(x)=-61$, there are members with no assigned days and only 19\% of members have 2 days scheduled. In the other solution, $z_4(x)=-125$, a higher percentage of members have 2 days assigned, 28\%, and 8\% of them have 3 or more. This solution does not have members without any days assigned.

\subsection{Results}\label{res}

\noindent Since the initial solutions are already generated \textit{via} crossovers between good quality solutions and there are a lot of interdependencies and relationships between certain assignments, the genetic algorithm for this case study runs for only 5 generations. The population size is 136 individuals and instead of using only the partial solutions from the last generation, the solutions from all generations are considered in the second stage. The time limit for all optimisation iterations is 120 seconds. 

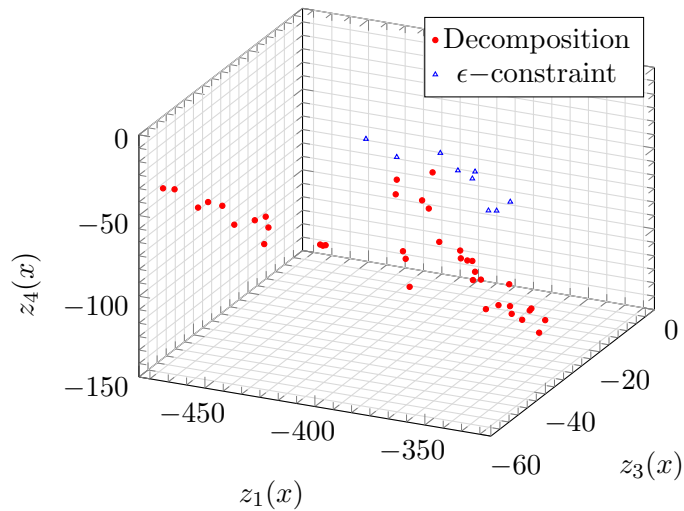
\begin{figure}[H]
    \centering
    \begin{tikzpicture}
    \begin{axis}[
        xlabel={$z_1(x)$},
        ylabel={$z_3(x)$},
        zlabel={$z_4(x)$},
        xmin=-480, xmax=-320,
        ymin=-60, ymax=0,
        zmin=-150, zmax=0,
        grid=both, 
        grid style={gray!30}, 
        minor tick num=5,
        ztick distance = 50,
    ]
    
    \addplot3[
        only marks,
        mark=*,
        mark size=1pt,
        red
    ] coordinates {
     
(-426.0,-9.0,-82.0)
(-474.0,-56.0,-37.0)
(-330.0,-32.0,-117.0)
(-452.0,-52.0,-48.0)
(-470.0,-55.0,-38.0)
(-468.0,-48.0,-58.0)
(-404.0,-15.0,-87.0)
(-340.0,-29.0,-116.0)
(-334.0,-31.0,-127.0)
(-342.0,-28.0,-119.0)
(-360.0,-21.0,-116.0)
(-378.0,-20.0,-107.0)
(-378.0,-19.0,-115.0)
(-386.0,-18.0,-105.0)
(-374.0,-23.0,-114.0)
(-382.0,-21.0,-105.0)
(-384.0,-17.0,-112.0)
(-398.0,-16.0,-105.0)
(-412.0,-11.0,-89.0)
(-406.0,-12.0,-69.0)
(-424.0,-11.0,-88.0)
(-454.0,-46.0,-68.0)
(-428.0,-56.0,-61.0)
(-436.0,-49.0,-55.0)
(-436.0,-48.0,-63.0)
(-346.0,-32.0,-112.0)
(-436.0,-53.0,-52.0)
(-350.0,-33.0,-111.0)
(-456.0,-54.0,-44.0)
(-386.0,-39.0,-78.0)
(-386.0,-38.0,-84.0)
(-352.0,-36.0,-110.0)
(-358.0,-33.0,-97.0)
(-410.0,-50.0,-65.0)
(-410.0,-49.0,-67.0)
(-410.0,-48.0,-68.0)
(-338.0,-34.0,-116.0)
(-378.0,-43.0,-93.0)
(-344.0,-33.0,-115.0)

    };
     \addplot3[only marks,
        mark=triangle,
        mark size=1pt,
        blue
    ] coordinates {
(	-378	,	-20,	-56		)
(	-378		,	-19,	-53	)
(	-362		,	-19,	-68	)
(	-400	,	-14,	-53		)
(	-440		,	-9,	-60	)
(	-426	,	-9	,	-68	)
(	-392		,	-14,	-62	)
(	-362		,	-27,	-63	)
(	-362		,	-24,	-67	)

    };
    \legend{Decomposition, $\epsilon-$constraint}
    \end{axis}
\end{tikzpicture}
    \caption{Non-dominated solutions for the DEG case study}
    \label{nds3d}
\end{figure}

The total runtime for the decomposition experiment was 1600 seconds and for the $\epsilon$-constraint it was 2400 seconds. The decomposition experiment found 39 non-dominated solutions and the augmented $\epsilon$-constraint found 9. These non-dominated solutions are presented in Figure \ref{nds3d}. 

\begin{figure}[H]
    \centering
    
    \begin{subfigure}[b]{0.45\textwidth}
        \centering
       \scalebox{0.8}{ \begin{tikzpicture}
                \begin{axis}[
        xlabel={$z_1(x)$},
        ylabel={$z_3(x)$},
        xmin=-480, xmax=-320,
        ymin=-60, ymax=0,
        xtick distance=40,
        grid=both, 
        grid style={gray!30}, 
        minor tick num=1,
    ]
    
    \addplot[
        only marks,
        mark=*,
        mark size=1pt,
        red
    ] coordinates {
     
(-426.0,-9.0)
(-474.0,-56.0)
(-330.0,-32.0)
(-452.0,-52.0)
(-470.0,-55.0)
(-468.0,-48.0)
(-404.0,-15.0)
(-340.0,-29.0)
(-334.0,-31.0)
(-342.0,-28.0)
(-360.0,-21.0)
(-378.0,-20.0)
(-378.0,-19.0)
(-386.0,-18.0)
(-374.0,-23.0)
(-382.0,-21.0)
(-384.0,-17.0)
(-398.0,-16.0)
(-412.0,-11.0)
(-406.0,-12.0)
(-424.0,-11.0)
(-454.0,-46.0)
(-428.0,-56.0)
(-436.0,-49.0)
(-436.0,-48.0)
(-346.0,-32.0)
(-436.0,-53.0)
(-350.0,-33.0)
(-456.0,-54.0)
(-386.0,-39.0)
(-386.0,-38.0)
(-352.0,-36.0)
(-358.0,-33.0)
(-410.0,-50.0)
(-410.0,-49.0)
(-410.0,-48.0)
(-338.0,-34.0)
(-378.0,-43.0)
(-344.0,-33.0)
            };

                \addplot[only marks,
        mark=triangle,
        mark size=1pt,
        blue
    ] coordinates {
(	-378	,	-20	)
(	-378		,	-19	)
(	-362		,	-19	)
(	-400	,	-14		)
(	-440		,	-9	)
(	-426	,	-9		)
(	-392		,	-14	)
(	-362		,	-27	)
(	-362		,	-24	)

    };
            \end{axis}
        \end{tikzpicture}}
        \caption{Pairwise trade-off profile between $z_1(x)$ and $z_3(x)$}
        \label{fig:sub1}
    \end{subfigure}
    \hfill
    \begin{subfigure}[b]{0.45\textwidth}
        \centering
        \scalebox{0.8}{\begin{tikzpicture}
                \begin{axis}[
        xlabel={$z_1(x)$},
        ylabel={$z_4(x)$},
        xmin=-480, xmax=-320,
        ymin=-150, ymax=0,
       xtick distance=40,
        grid=both, 
        grid style={gray!30}, 
        minor tick num=1,
    ]
    
    \addplot[
        only marks,
        mark=*,
        mark size=1pt,
        red
    ] coordinates {
     
(-426.0,-82.0)
(-474.0,-37.0)
(-330.0,-117.0)
(-452.0,-48.0)
(-470.0,-38.0)
(-468.0,-58.0)
(-404.0,-87.0)
(-340.0,-116.0)
(-334.0,-127.0)
(-342.0,-119.0)
(-360.0,-116.0)
(-378.0,-107.0)
(-378.0,-115.0)
(-386.0,-105.0)
(-374.0,-114.0)
(-382.0,-105.0)
(-384.0,-112.0)
(-398.0,-105.0)
(-412.0,-89.0)
(-406.0,-69.0)
(-424.0,-88.0)
(-454.0,-68.0)
(-428.0,-61.0)
(-436.0,-55.0)
(-436.0,-63.0)
(-346.0,-112.0)
(-436.0,-52.0)
(-350.0,-111.0)
(-456.0,-44.0)
(-386.0,-78.0)
(-386.0,-84.0)
(-352.0,-110.0)
(-358.0,-97.0)
(-410.0,-65.0)
(-410.0,-67.0)
(-410.0,-68.0)
(-338.0,-116.0)
(-378.0,-93.0)
(-344.0,-115.0)

    };
     \addplot[only marks,
        mark=triangle,
        mark size=1pt,
        blue
    ] coordinates {
(	-378,	-56		)
(	-378,	-53	)
(	-362,	-68	)
(	-400,	-53		)
(	-440,	-60	)
(	-426,	-68	)
(	-392,	-62	)
(	-362,	-63	)
(	-362,	-67	)

    };
            \end{axis}
        \end{tikzpicture}}
        \caption{Pairwise trade-off profile between $z_1(x)$ and $z_4(x)$}
        \label{fig:sub2}
    \end{subfigure}
    \hfill
    \begin{subfigure}[b]{0.5\textwidth}
        \centering
        \scalebox{0.8}{\begin{tikzpicture}
            \begin{axis}[
        xlabel={$z_3(x)$},
        ylabel={$z_4(x)$},
        xmin=-60, xmax=0,
        ymin=-150, ymax=0,
        grid=both, 
        grid style={gray!30}, 
        minor tick num=1,
    ]
    
    \addplot[
        only marks,
        mark=*,
        mark size=1pt,
        red
    ] coordinates {
    
(-9.0,-82.0)
(-56.0,-37.0)
(-32.0,-117.0)
(-52.0,-48.0)
(-55.0,-38.0)
(-48.0,-58.0)
(-15.0,-87.0)
(-29.0,-116.0)
(-31.0,-127.0)
(-28.0,-119.0)
(-21.0,-116.0)
(-20.0,-107.0)
(-19.0,-115.0)
(-18.0,-105.0)
(-23.0,-114.0)
(-21.0,-105.0)
(-17.0,-112.0)
(-16.0,-105.0)
(-11.0,-89.0)
(-12.0,-69.0)
(-11.0,-88.0)
(-46.0,-68.0)
(-56.0,-61.0)
(-49.0,-55.0)
(-48.0,-63.0)
(-32.0,-112.0)
(-53.0,-52.0)
(-33.0,-111.0)
(-54.0,-44.0)
(-39.0,-78.0)
(-38.0,-84.0)
(-36.0,-110.0)
(-33.0,-97.0)
(-50.0,-65.0)
(-49.0,-67.0)
(-48.0,-68.0)
(-34.0,-116.0)
(-43.0,-93.0)
(-33.0,-115.0)

    };
     \addplot[only marks,
        mark=triangle,
        mark size=1pt,
        blue
    ] coordinates {
(	-20,	-56		)
(	-19,	-53	)
(	-19,	-68	)
(	-14,	-53		)
(	-9,	-60	)
(	-9	,	-68	)
(	-14,	-62	)
(	-27,	-63	)
(	-24,	-67	)

    };
    \legend{Decomposition, $\epsilon-$constraint}
            
            \end{axis}
        \end{tikzpicture}}
        \caption{Pairwise trade-off profile between $z_3(x)$ and $z_4(x)$}
        \label{fig:sub3}
    \end{subfigure}
    
    \caption{Pairwise trade-off profiles between the three objective functions}
    \label{fig:subfigures}
\end{figure}
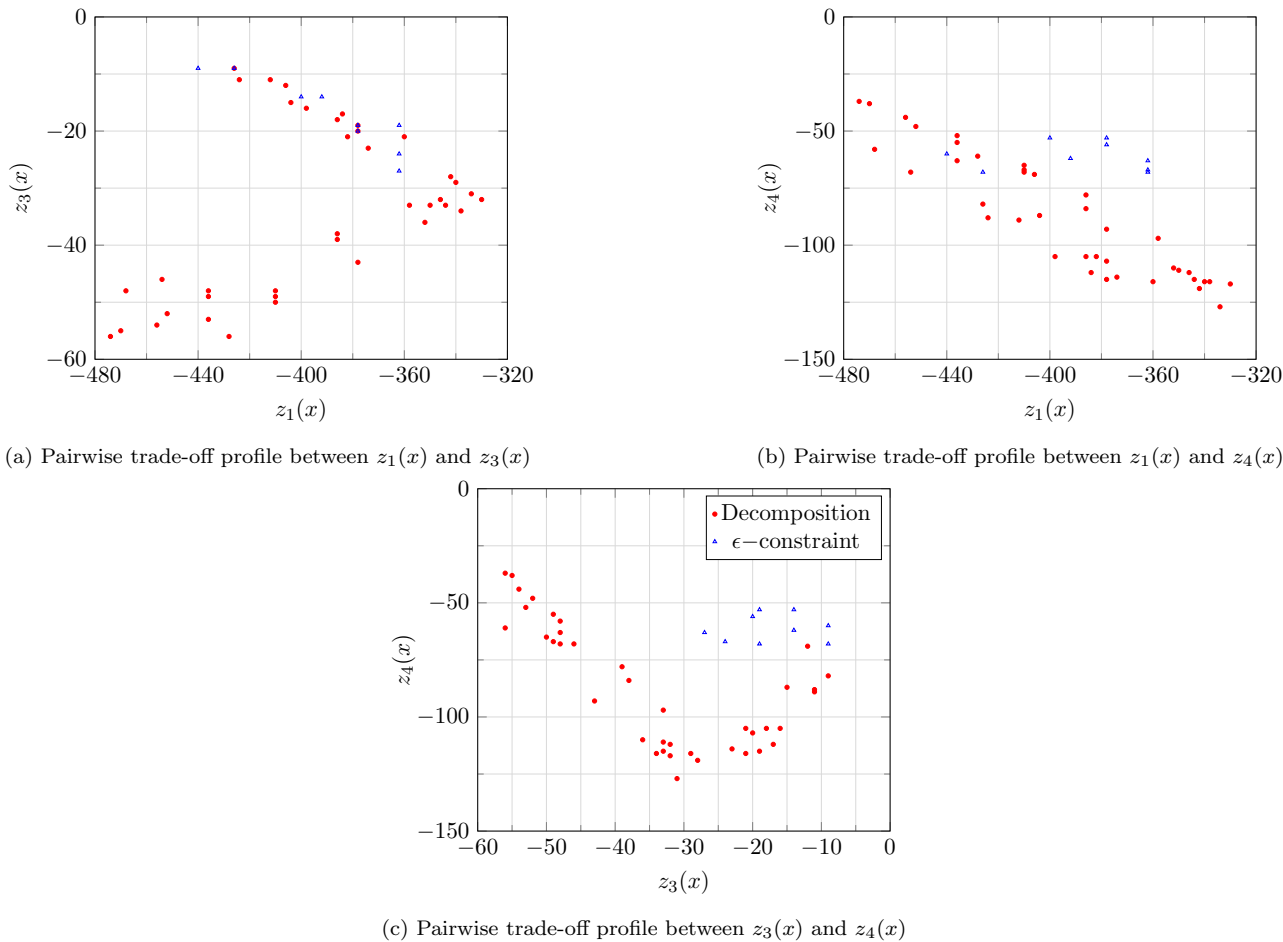

For a clearer visual representation, pairwise trade-off profiles between the three objectives are presented in Figure \ref{fig:subfigures}.

Looking separately at each objective, our method found a solution with better performance than the $\epsilon$-constraint method (the maximum value for $z_3(x)$ with the $\epsilon$-constraint is -9 and with decomposition it is -8). Let us note that this happens because the time limit \textit{per} iteration stops the $\epsilon$-constraint iterations before the gap reaches 0. 

The workload balance objective, $z_1(x)$, is the objective where our method performed the best. This is expected as it is the first stage objective which is maximised by the genetic algorithm. Regarding the other two objectives, the time-slot preference objective, $z_3(x)$, appears to have a better performance than the committee days objective, $z_4(x)$. Accordingly, the non-dominated solutions considering the pairing between the first stage objective, $z_1(x)$, and the time-slot preference objective, $z_3(x)$, share very similar trade-off profiles for the two methods. However, when considering pairings with the committee days objective, $z_4(x)$, several non-dominated solutions found by our method would be dominated by those found by the $\epsilon$-constraint method.

In the first stage, the objectives being regarded are the workload balance objective, $z_1(x)$, and a proxy objective that maximises the number of time slots where the committee members being assigned to the defence are all available. Then in the second stage, the time-slot preference objective, $z_3(x)$, and the committee days objective, $z_4(x)$, are considered. The differences in performance of these last objectives in the case study can be explained by how well the proxy objective predicts their performance. For the time-slot preference objective, $z_3(x)$, having more time slots available for a defence strongly increases the likelihood of a slot where all members prefer to have a defence scheduled. However, for the committee days objective, $z_4(x)$, this relationship is not as strong. For example, a member might be present in two committees that have many available time slots, but few of these slots may be on the same day. This means that the member would have to be scheduled for two different days. 

In this case study, one of the monolithic objectives is assessed in the first stage and two are assessed in the second stage. The performance of the method should improve when more monolithic objectives are considered in the first stage. Nonetheless, for one of the other objectives the performance is also very good, as it has a strong relation with the proxy objective. This proves that the method can yield excellent results if good quality proxy objectives can be defined. Perhaps, if another proxy objective with a stronger relationship with the committee days objective is added in the first stage, its performance could be improved to the same level as the time slot preference objective.

To sum up, we can derive the following conclusions from this case study:
\begin{itemize}[label={--}]
\item Our method improved the upper bounds of each objective function;
\item First stage objectives can be expected to have better performances than second stage objectives;
\item The performance of solutions regarding the second stage objective is affected by how strong the relationship they have with the proxy objectives.
\end{itemize}
\section{Conclusion}\label{conc}

\noindent This work proposes a decomposed multi-objective method combining genetic algorithms and an augmented $\epsilon$-constraint method. It involves breaking down the monolithic problem into a series of sequential multi-objective problems. Each partial solution obtained in the preceding stage is subjected to solving a multi-objective problem, enhancing the overall efficiency of the method. Furthermore, proxy objective functions are incorporated in earlier stages to forecast the performance of objectives in subsequent stages.

This method can be applied to problems with specific characteristics. They must be multi-objective problems, must be decomposable into sequential sub-problems, and the partial solutions of previous sub-problems must be able to be translated into proxy objectives for the performance of the next stages. The thesis defence scheduling problem fits these criteria. In our application, the first stage is a genetic algorithm which finds several committee configurations. In the second stage, an augmented $\epsilon$-constraint method is ran for several configurations found in the previous stage.

The method undergoes testing across both small and large size randomly generated instances, as well as a real-world case study. In smaller size instances, where the augmented $\epsilon$-constraint iterations can achieve optimality within a 120-second limit, the decomposed approach consumes between 8\% and 32\% of the time monolithic resolutions do. However, monolithic resolutions excel in delivering higher-quality solutions, which is anticipated given the comparison between a partial optimization technique and a deterministic method capable of reaching optimality. For larger size instances where optimality is not achieved by the augmented $\epsilon$-constraint iterations, the decomposed method requires between 6\% and 18\% of the time. Furthermore, it demonstrates an ability to identify non-dominated sets with greater hyper-volume indicator values.
The real-world case study highlights the performance across three objectives: one evaluated in the initial stage and two in the subsequent stage. The performance of the objective in the initial stage, along with one of the objectives in the subsequent stage, closely aligns with monolithic resolutions. However, the other objective in the subsequent stage exhibits slightly inferior performance. This discrepancy is attributed to one of the objectives in the subsequent stage having a stronger correlation with the proxy objective compared to the other.

This is an efficient method which can produce excellent results for multi-objective problems. Regarding future research, it should be interesting to see applications with different problems and using different multi-objective algorithms. Moreover, in our tests, a single proxy objective was used to predict the performance of the next stage. Nonetheless, approaches using multiple different proxy objectives focused on each objective might yield better results.

\section*{Acknowledgements}
\noindent The corresponding author acknowledges the Portuguese national funds through the FCT - Foundation for Science and Technology, I.P., under the project UI/BD/154023/2022. All authors acknowledge the support by national funds through FCT (Fundação para a Ciência e a Tecnologia), Portugal under the project UIDB/00097/2020. Declarations of interest: none. 
\bibliographystyle{model2-names}

\end{document}